\documentclass[11pt]{article}

\usepackage{fancyhdr}
\usepackage{epsfig}
\usepackage{color}

\usepackage{graphicx}
\usepackage{latexsym}
\usepackage{amssymb}
\usepackage{amsfonts}
\usepackage{amsmath}
\usepackage{amsthm}
\usepackage{paralist}
\usepackage{bm}
\usepackage{caption}
\usepackage{subcaption}
\usepackage[font=small,skip=0pt]{caption}
\usepackage{float}

\theoremstyle{plain}
\newtheorem{res}{Result}
\theoremstyle{definition}
\newtheorem{rem}{Remark}

\DeclareMathOperator{\PP}{{\bf P}}
\DeclareMathOperator{\E}{{\bf E}}
\DeclareMathOperator{\Var}{{\bf Var}}
\DeclareMathOperator{\Cov}{{\bf Cov}}
\DeclareMathOperator{\PPe}{{\bf P}_{\eta}}
\DeclareMathOperator{\Ee}{{\bf E}_{\eta}}
\DeclareMathOperator{\Vare}{{\bf Var}_{\eta}}

\DeclareMathOperator{\Ete}{{\bf E}_{\tilde\eta}}
\DeclareMathOperator{\Varte}{{\bf Var}_{\tilde\eta}}

\DeclareMathOperator{\Ex}{{\bf E}_{\xi}}
\DeclareMathOperator{\Varx}{{\bf Var}_{\xi}}

\newcommand{\cf}{\mbox{cf.\frenchspacing}}

\usepackage[top=1in, bottom=1in, left=1in, right=1in]{geometry}
\begin{document}

\LARGE
\begin{center}
\textbf {A solution in small area estimation problems}\\[2.5\baselineskip]

\small
\text{Andrius {\v C}iginas and Tomas Rudys}\\

\medskip

{\footnotesize
Vilnius University Institute of Mathematics and Informatics, LT-08663 Vilnius, Lithuania \\[2.5\baselineskip]}
\end{center}

\small
\begin{abstract}

We present a new method in problems where estimates are needed for finite population domains with small or even zero sample sizes. In contrast to known estimation methods, an auxiliary information is used to model sizes of population units instead of a direct prediction of their values of interest. In particular, via an additional characterization of regression models, we incorporate a scatter and variabilities of the units sizes into an estimator, and then it uses an information of the whole sample by taking into an account a location of the estimation domain inside the population. To reduce an impact of the introduced domain total estimator bias to the mean square error, we construct also a regression type version of the estimator. An efficiency of the method proposed  is shown in a simulation study.
   
\end{abstract}
\vskip 4mm

\normalsize

\noindent\textbf{Keywords:} small area estimation, auxiliary information, linear regression, order statistics.

\vskip 2mm

\noindent\textbf{MSC classes:} 62D05 (62J05)

\section{Introduction}\label{s:1}

\let\thefootnote\relax\footnotetext{The research of the first author is partially supported by European Union Structural Funds project "Postdoctoral Fellowship Implementation in Lithuania".} 

The topic small area estimation (SAE) becomes more and more actual and thus popular in the last decades because of a need to get more inferences about survey populations than it seems possible with given samples or even after a careful planning of them. Actually, historically established but misleading term `SAE' means an estimation in the finite population domain (not necessarily a small) where sample size is not sufficient to get estimates of an adequate precision. More specifically, problems of this kind start to arise when we need to plan the sample in order to estimate population characteristics in the number of the population subsets which are determined by various classifications of the population and which are intersecting in various ways. A sampling design, ensuring a sample in each of the intersections, leads to the large total sample size and thus increases a cost of the survey. Next, if samples obtained are small in domains, applications of the classical estimation theory usually fail, in the sense of estimates quality, and such a failure not so much depends on a good auxiliary information availability. In SAE, auxiliary information plays a crucial role similarly as it is important in the traditional survey sampling where it is extensively used via ratio, regression, calibration estimators, etc., see, e.g., \cite{C_1977, SSW_1992, DS_1992, F_2009}. The main difference is that, in distinctive SAE methods, the auxiliary data are used to extract an information about the variable or parameter of interest from the sample outside the domain (so-called indirect estimation), while, in the classical case, the sample information of the domain of interest is incorporated into estimation only (direct estimation).

SAE methods can be classified roughly into two large groups: design-based and model-based methods. The design-based methodologies correspond to an understanding of the classical estimation. Data models, included into an estimation of this type, are advisedly used to keep the basic properties of the resulting estimators such as an asymptotic consistency and  at least an approximate unbiasedness with respect to the sampling design, i.e., here a randomness is considered by the distribution of different samples appearance. A validity of the underlying model is not so much important, i.e., if the model is an incorrect then the variance of an initial estimator can be left unreduced, but, in any case, the consistency and an acceptable unbiasedness of the estimator still hold. Here the underlying model is called therefore assisting. See \cite{LV_2009} and also \cite{P_2013}, for a wide review of the design-based methods. However, the unbiasedness of the design-based estimators has its cost, i.e., for the population domains with markedly small sample size, the variances of estimators are often large. Then estimators from the class of model-based estimators can be a better choice, as it is pointed, e.g., in \cite{R_2005}. The number of model-based methods are collected in \cite{R_2003}, see also \cite{P_2013} and \cite{D_2009} on recent results. We refer also to reviews in \cite{GR_1994} and \cite{P_2002}.  Differently from the design-based methods, the mean square error (MSE) of a model-based estimator is defined and estimated with a respect to the model of an explicit form.

To explain better a place of our estimation approach, let us discuss shortly about two popular types of small area estimators: synthetic and empirical best linear unbiased predictor (EBLUP). The former could be placed between the classes described above, because its properties can be considered from the point of the sampling design, while it is indirect, but uses data models implicitly. The synthetic estimator works well if the population (or a larger area, containing the domain of interest) has similar characteristics of a model to that of the domain. Estimators of the EBLUP type are the members of the model-based estimators class. They take into an account a particularity of the domain. Estimators, we propose in the paper, have the mentioned properties of synthetic estimators except that effects of the estimation domain are incorporated by a new way which is in a sense similar to that in the EBLUP case. 

Let us turn to specific assumptions. Consider a finite population ${\cal U}=\{ 1,\ldots, N\}$ of size $N$. Assume that, in order to estimate a parameter of the population, the sample $s=\{i_1,\ldots, i_n\}$ of size $n$ is drawn from ${\cal U}$, according to a sampling design $p(\cdot)$. Here $p(s)$ is the probability to get the particular $s$. Let $\pi_i=\PP\{i\in s\}>0$ and $\pi_{ij}=\PP\{i,j\in s\}>0$ be the inclusion into the sample probabilities for the population element $i$ and for the pair of elements $i$ and $j$. Here the operator $\PP$ and further $\E$, $\Var$ and $\Cov$ mean the probability, expectation, variance and covariance according to the sampling design, respectively. It is a common situation when a construction of the design $p(\cdot)$ is related closely to an auxiliary variable, say, $z$ with values $\{z_1,\ldots, z_N \}$ known for all units of ${\cal U}$. In particular, e.g., for a stratified simple random sampling design and for probability proportional-to-size sampling, the inclusion probability $\pi_i$ represents an importance of the population unit $i$ by the relative size of $z_i$. We renumber the population ${\cal U}$ in order to have $z_1\leq \cdots \leq z_N$. Let $y$ be the variable of interest with the fixed values $\{y_1,\ldots, y_N\}$ in the population, and, we aim to estimate the total 
\begin{equation}\label{sum_yD}
t_{y;{\cal D}}=\sum_{i\in {\cal D}} y_i
\end{equation}
where ${\cal D}\subseteq {\cal U}$ is any non-empty set. If the particular estimation domain ${\cal D}$ is known before the sample selection, then the sampling design  $p(\cdot)$ can be organized to get a sufficient for quality requirements of estimates sample size in that domain. But if we are interested in  ${\cal D}$ after the sample selection and collection of the data $y_{i_1},\ldots, y_{i_n}$, the sample size can be too small in ${\cal D}$, and this leads to get bad quality estimates, if the estimators of \eqref{sum_yD} are, e.g., the direct Horvitz--Thompson (H--T) estimator
\begin{equation}\label{HT}
\hat t_{y;{\cal D}}^{HT}=\sum_{i\in s\cap {\cal D}}d_i y_i, \quad \text{where} \quad  d_i=\frac{1}{\pi_i},
\end{equation}
or a direct generalized regression (GREG) estimator which can be also represented by a similar to \eqref{HT} form. Furthermore, in the case of empty sample, the direct estimators cannot be applied at all. Assume that at the estimation stage, we have a different from $z$ auxiliary variable $x$ and all its values $\{x_1,\ldots, x_N \}$ are known. Let these values be the realization of independent random variables $X_1,\ldots, X_N$ modeled by the linear regression model 
\begin{equation}\label{x_reg_z}
X_i=\alpha_1+\alpha_2 z_i+\delta_i, \qquad i=1,\ldots, N,
\end{equation}
which we call $\eta$. Here, in general, $X_i$ has the distribution function $F_i(\cdot)$ with $\Ee X_i=\alpha_1+\alpha_2 z_i$ and $\Vare X_i=\tau_i^2$, where the symbols $\Ee$ and $\Vare$ denote expectation and variance with respect to the model $\eta$. Thus we assume a superpopulation model, where $\alpha_1$, $\alpha_2$ and $\tau_1,\ldots, \tau_N$ are unknown model parameters.  In traditional survey sampling estimation problems, an introduction of the auxiliary variable $x$ means usually that it is a better than $z$ linear predictor of $y$. Then the superpopulation model
\begin{equation}\label{y_reg_x}
Y_i=\beta_1+\beta_2 x_i+\varepsilon_i, \qquad i=1,\ldots, N,
\end{equation}
we call $\xi$, is used in an estimation. Here $Y_1,\ldots, Y_N$ are independent random versions of the fixed $y$ values with the expectations $\Ex Y_i=\beta_1+\beta_2 x_i$ and variances  $\Varx Y_i=\sigma_i^2$. To be consecutive, we do not treat model \eqref{y_reg_x} as an important in Section \ref{s:2}, where we keep in mind only that both $z$ and $x$ more or less represent a size of $y$. Next, in Section \ref{s:3}, we take into an account linear relation \eqref{y_reg_x}.

Estimators of the synthetic and EBLUP types are based on a common principle to predict $y$ values in the domain of interest using linear relations between $y$ and auxiliary variables by the data from `neighbour' population areas. The approach, proposed in Section \ref{s:2}, is different. Here we model sizes of the population elements instead of their $y$ values. More specifically, we introduce an estimator of the form
\begin{equation*}
\hat t_{y;{\cal D}}=\sum_{i\in s} \widehat w_iy_i, \quad \text{where} \quad  \widehat w_i=\hat\theta_{i;{\cal D}}d_i, \quad 0\leq\hat\theta_{i;{\cal D}}\leq 1,
\end{equation*}
similar to those in the design-based estimation, \cf{} \eqref{HT}. Here the numbers $\hat\theta_{i;{\cal D}}$, $i\in s$, we get from model \eqref{x_reg_z}, describe how the sample elements are similar to elements of the domain ${\cal D}$. Such an incorporation of the sizes model \eqref{x_reg_z} into the estimation can be interpreted, for instance, as a mimic of a change of the population over time or as a treatment of a size relativity. Therefore,  we give the name \emph{hidden randomness} (HR) to the method, which we present detailed in Section \ref{s:2}. Further, in Section \ref{s:3}, we aim to reduce the bias of the HR estimator by constructing a regression type version of the estimator. In Section \ref{s:4}, we present the simulation study, where we compare the introduced estimators with the synthetic, EBLUP and direct GREG estimators. Conclusions are given in Section \ref{s:5}, where a relaxation of auxiliary model (assumption) \eqref{x_reg_z} is discussed additionally.

\section{Method of hidden randomness and small area estimator}\label{s:2}

Let $X_{(1)}\leq \cdots \leq X_{(N)}$ denote the order statistics of random variables $X_1,\ldots, X_N$ generated by linear regression model \eqref{x_reg_z}. We introduce the probabilities
\begin{equation}\label{p_ij}
p_{ij}=\PPe\left\{ X_i=X_{(j)}\right\}, \qquad i,j=1,\ldots, N.  
\end{equation}
These numbers are the parameters of the model $\eta$. They depend on the basic parameters $\alpha_2$ and $\tau_1,\ldots, \tau_N$, but characterize also the model additionally. In particular, characteristics \eqref{p_ij} contain an information about the scatter of the given $z$ values in ${\cal U}$. Here, for the fixed unit $i\in {\cal U}$, we interpret the probabilities in the following way: the collection $\{p_{ij}, j=1,\ldots, N\}$ shows how the initial size $z_i$ tends to variate. Define the numbers 
\begin{equation*}\label{CP}
t_j=\sum_{i\in {\cal U}} p_{ij}y_i, \qquad j=1,\ldots, N.
\end{equation*}
Let ${\cal D}=\{ j_1,\ldots, j_m\}\subseteq {\cal U}$ be the domain of interest of size $1\leq m\leq N$. Since the identity $\sum_{i\in {\cal U}} y_i=\sum_{j\in {\cal U}} t_j$ is satisfied, we approximate target sum \eqref{sum_yD} as follows:
\begin{equation*}\label{sum_appr}
t_{y;{\cal D}}\approx \sum_{j\in {\cal D}} t_j=\sum_{i\in {\cal U}}\theta_{i;{\cal D}}y_i, \quad \text{where} \quad  \theta_{i;{\cal D}}=\sum_{j\in {\cal D}}p_{ij}
\end{equation*}
can be treated as a probability that the population element $i$ should `belong' to the domain ${\cal D}$. We define the HR estimator of the sum $t_{y;{\cal D}}$ by
\begin{equation}\label{HT_type}
\hat t_{y;{\cal D}}^{HR}=\sum_{i\in s}\hat\theta_{i;{\cal D}}d_iy_i, \quad \text{with} \quad  \hat\theta_{i;{\cal D}}=\sum_{j\in {\cal D}}\hat p_{ij}.
\end{equation}
Here estimates $\hat p_{ij}$ are plugged-in instead of the probabilities $p_{ij}$, because the parameters $\alpha_1$, $\alpha_2$ and $\tau_1,\ldots, \tau_N$ of model \eqref{x_reg_z} are assumed not known as well as $F_1(\cdot), \ldots, F_N(\cdot)$ are not specified. An exact calculation of \eqref{p_ij} is complicated even the mentioned model characteristics are known, since, except trivial cases, the random variables $X_1,\ldots, X_N$ are non-identically distributed. For instance, evaluations of distributions of the corresponding order statistics require intensive computing, see, e.g., \cite{CHB_2009} and references therein. Therefore, we propose two alternative ways to evaluate the numbers $\theta_{i;{\cal D}}$, $i=1,\ldots, N$, see Appendix \ref{a:1}.
\begin{rem}\label{rem:1}
In the case of ${\cal D}={\cal U}$, estimator \eqref{HT_type} coincides with the design unbiased H--T estimator. Two other connections with known estimators are found in the following separate cases of  \eqref{HT_type}. First, assuming that the ranks of $X_1, \ldots, X_N$ satisfy $\{R_1,\ldots, R_N\}\equiv\{1,\ldots, N\}$ or that the coefficient of correlation between the variables $x$ and $z$ is $\rho_{xz}=1$, we get $p_{ij}=\mathbb{I}\{ i=j\}$, where $ \mathbb{I}\{ \cdot\}$ is the indicator function. Then \eqref{HT_type} yields direct H--T estimator \eqref{HT}. Second, taking $\alpha_2=0$, $\tau_i=\tau>0$ and $F_i(\cdot)\equiv F(\cdot)$ in model \eqref{x_reg_z}, we obtain $p_{ij}=N^{-1}$. With this absence of the linear dependence between $x$ and $z$, we arrive from \eqref{HT_type} to
\begin{equation}\label{simp_syn}
\hat t_{y;{\cal D}}^{S}=\frac{m}{N}\sum_{i\in s} d_i y_i,
\end{equation}
which is the simplest synthetic estimator, where the additional information $x$ is not used (or is not helpful) at the estimation. In comparison with \eqref{HT}, estimator \eqref{simp_syn} has small variance, but its bias can be large if the domain ${\cal D}$  is not homogeneous with respect to the population. 
\end{rem}

Let us consider the MSE of estimator \eqref{HT_type}. The indirect estimators are generally design biased, and the construction of \eqref{HT_type} implies a biasedness too. While the variance part of the MSE is easily estimated using standard design-based methods, an estimation of the bias is more difficult. There are specific MSE estimation methods for synthetic estimators, see, e.g., \cite{R_2003}, which can be applied also to the HR estimator, but here, we follow a course of our methodology. The bias of \eqref{HT_type} has the expression
\begin{equation}\label{HT_bias}
B_y^{HR}=\mathrm{BIAS}(\hat t_{y;{\cal D}}^{HR})=\sum_{i\in {\cal U}}\hat\theta_{i;{\cal D}} y_i - \sum_{i\in {\cal D}} y_i
\end{equation}
depending on parameter \eqref{sum_yD}. Therefore, the design unbiased estimator of \eqref{HT_bias} is not suitable. We introduce the following:
\begin{equation}\label{bias_est}
\widehat B_{y}^{HR}=\frac{1-\rho_{xz}^2}{\rho_{xz}^2} \sum_{i\in s}\left(\frac{m}{N}-\hat\theta_{i;{\cal D}}\right)d_i y_i.
\end{equation}
This estimator is consistent with the first two special cases of Remark \ref{rem:1}, where bias \eqref{HT_bias} equals zero. In the case of $\rho_{xz}=0$ only, estimator \eqref{bias_est} is not clearly defined, because of a complexity of (properties of) parameters \eqref{p_ij}. To summarize, we formulate the statement on the accuracy of estimator \eqref{HT_type}. For short, we denote $a_{ij}=d_id_j\pi_{ij}-1$ with $a_{ii}=d_i-1$. Similarly, we write $\tilde a_{ij}=a_{ij}/\pi_{ij}$.  
\begin{res}\label{res:1}
\emph{(i)} The \emph{MSE} of the estimator $\hat t_{y;{\cal D}}^{HR}$ of the sum $t_{y;{\cal D}}$ in the domain ${\cal D}$ is 
\begin{equation}\label{HT_mse}
\mathrm{MSE}(\hat t_{y;{\cal D}}^{HR})=\Var \hat t_{y;{\cal D}}^{HR} + \bigl( B_y^{HR}\bigr)^2,
\end{equation}
where $\Var \hat t_{y;{\cal D}}^{HR}=\sum_{i,j\in {\cal U}} \hat\theta_{i;{\cal D}} \hat\theta_{j;{\cal D}} a_{ij}y_iy_j$, and $B_y^{HR}$ is given by \eqref{HT_bias}.
\smallskip

\emph{(ii)} The estimator of \eqref{HT_mse} is
\begin{equation}\label{HT_mse_est}
\widehat{\mathrm{MSE}}(\hat t_{y;{\cal D}}^{HR})=\widehat{\Var} \hat t_{y;{\cal D}}^{HR} + \bigl(\widehat B_{y}^{HR} \bigr)^2,
\end{equation}
where $\widehat{\Var} \hat t_{y;{\cal D}}^{HR}=\sum_{i,j\in s} \hat\theta_{i;{\cal D}} \hat\theta_{j;{\cal D}} \tilde a_{ij}y_iy_j$ is the unbiased estimator of $\Var \hat t_{y;{\cal D}}^{HR}$, and $\widehat B_{y}^{HR}$ is in \eqref{bias_est}.
\end{res}
\begin{rem}
If the total sample size $n$ and the size $m$ of the domain ${\cal D}$ are small, then the random variable $\hat t_{1;{\cal D}}^{HR}:=\sum_{i\in s} \hat\theta_{i;{\cal D}} d_i$, induced by the sampling design, varies more about $m$. Therefore, applying the ratio estimator 
\begin{equation*}\label{HT1_type}
\hat t_{y;{\cal D}}^{HR1}=m\frac{\hat t_{y;{\cal D}}^{HR}}{\hat t_{1;{\cal D}}^{HR}},
\end{equation*}
one can expect to reduce an impact of this error source to the MSE of the HR estimator. 
\end{rem}

\section{Regression type estimator}\label{s:3}

Assume that the auxiliary variable $x$ explains the study variable $y$ well according to model \eqref{y_reg_x}. By Result \ref{res:1}, the variance of \eqref{HT_type} does not exceed the variance of the direct whole population total  H--T estimator. Therefore, we aim to exploit the variable $x$ in a reduction of the HR estimator bias by considering a regression type estimator of the form $\hat t_{y;{\cal D}}^{HR}+b(t_{x;{\cal D}}-\hat t_{x;{\cal D}}^{HR})$ with a properly chosen characteristic $b$. In particular, minimizing the MSE of this expression, we get 
\begin{equation}\label{b_0}
b_0=\frac{C_{yx}+B_y^{HR}B_x^{HR}}{V_x+\bigl(B_x^{HR}\bigr)^2}
\end{equation}
where $C_{yx}=\Cov(\hat t_{y;{\cal D}}^{HR}, \hat t_{x;{\cal D}}^{HR})=\sum_{i,j\in {\cal U}} \hat\theta_{i;{\cal D}} \hat\theta_{j;{\cal D}} a_{ij}y_ix_j$ and $V_x=\Var \hat t_{x;{\cal D}}^{HR}=C_{xx}$. 
Then we introduce the regression type HR (RHR) estimator
\begin{equation}\label{reg_type}
\hat t_{y;{\cal D}}^{RHR}=\hat t_{y;{\cal D}}^{HR}+\hat b_0\left(t_{x;{\cal D}}-\hat t_{x;{\cal D}}^{HR}\right),
\end{equation}
where 
\begin{equation}\label{b_0_est}
\hat b_0=\frac{\widehat C_{yx}+\widehat B_{y}^{HR} \widehat B_{x}^{HR}}{\widehat V_x+ \bigl(\widehat B_{x}^{HR}\bigr)^2}
\end{equation}
is the estimator of the parameter $b_0$. Here
\begin{equation*}
\widehat C_{yx}=\sum_{i,j\in s} \hat\theta_{i;{\cal D}} \hat\theta_{j;{\cal D}} \tilde a_{ij}y_ix_j \quad \text{and} \quad \widehat V_x=\sum_{i,j\in s} \hat\theta_{i;{\cal D}} \hat\theta_{j;{\cal D}} \tilde a_{ij}x_ix_j
\end{equation*}
are the design unbiased estimators of $C_{yx}$ and $V_x$, respectively.

To evaluate the MSE of estimator \eqref{reg_type}, we apply the traditional Taylor linearization to the estimators function $\hat t_{y;{\cal D}}^{RHR}$ at the point $(\E \hat t_{y;{\cal D}}^{HR}, \E \hat t_{x;{\cal D}}^{HR}, C_{yx}, V_x, B_y^{HR}, B_x^{HR})$ which yields
\begin{align*}
\begin{split}
\hat t_{y;{\cal D}}^{RHR}\approx \tilde t_{y;{\cal D}}^{RHR}:=&\hat t_{y;{\cal D}}^{HR}+b_0\bigl(t_{x;{\cal D}}-\hat t_{x;{\cal D}}^{HR}\bigr)-B_x^{HR}\left( V_x+\bigl( B_x^{HR}\bigr)^2\right)^{-1}\left\{ \widehat C_{yx}-C_{yx}-b_0\bigl( \widehat V_x-V_x\bigr)\right. \\
& \left. +B_x^{HR}\bigl( \widehat B_y^{HR}-B_y^{HR}\bigr)-\bigl( 2b_0B_x^{HR}-B_y^{HR}\bigr)\bigl( \widehat B_x^{HR}-B_x^{HR}\bigr)\right\}.
\end{split}
\end{align*}
According to this one-term Taylor approximation, the bias of estimator \eqref{reg_type} is approximated by the bias of $\tilde t_{y;{\cal D}}^{RHR}$. Next, one can verify that the latter bias consists of the term $B_y^{HR}-b_0 B_x^{HR}$ plus a remainder which is negligible if biases of the estimators $\widehat B_{y}^{HR}$ and $\widehat B_{x}^{HR}$ are of a smaller order than the corresponding biases $B_y^{HR}$ and $B_x^{HR}$. The same reasons imply that the variance of the term $\hat t_{y;{\cal D}}^{HR}+b_0(t_{x;{\cal D}}-\hat t_{x;{\cal D}}^{HR})$ dominates in the variance of $\tilde t_{y;{\cal D}}^{RHR}$. Therefore, we formulate the following result on the accuracy of estimator \eqref{reg_type}. 
\begin{res}\label{res:2}
\emph{(i)} The \emph{MSE} approximation of the estimator $\hat t_{y;{\cal D}}^{RHR}$ of the sum $t_{y;{\cal D}}$ in the domain ${\cal D}$ is 
\begin{equation}\label{RHR_mse}
\mathrm{MSE}(\hat t_{y;{\cal D}}^{RHR})\approx \Var \hat t_{y;{\cal D}}^{HR}+b_0^2 V_x-2b_0 C_{yx} + \bigl( B_y^{HR}-b_0B_x^{HR}\bigr)^2,
\end{equation}
where $\Var \hat t_{y;{\cal D}}^{HR}$ is the same as in \eqref{HT_mse}, $C_{yx}$ and $V_x$ are from expression \eqref{b_0} of $b_0$, and $B_y^{HR}$ and $B_x^{HR}$ are by formula \eqref{HT_bias}.
\smallskip

\emph{(ii)} The estimator of approximation \eqref{RHR_mse} is
\begin{equation}\label{RHR_mse_est}
\widehat{\mathrm{MSE}}(\hat t_{y;{\cal D}}^{RHR})=\widehat{\Var} \hat t_{y;{\cal D}}^{HR} + \hat b_0^2 \widehat V_x-2\hat b_0 \widehat C_{yx}+ \bigl(\widehat B_{y}^{HR}-\hat b_0 \widehat B_{x}^{HR} \bigr)^2,
\end{equation}
where $\widehat{\Var} \hat t_{y;{\cal D}}^{HR}$ is the same as in \eqref{HT_mse_est}, $\widehat C_{yx}$ and $\widehat V_x$ are from expression \eqref{b_0_est} of $\hat b_0$, and $\widehat B_{y}^{HR}$ and $\widehat B_{x}^{HR}$ are given by \eqref{bias_est}.
\end{res}
One can conclude from \eqref{RHR_mse} that the RHR estimator reduces bias \eqref{HT_bias} of estimator \eqref{HT_type} if the variables $y$ and $x$ are well-correlated.

\section{Simulation study}\label{s:4}

In this section, we compare RHR estimator \eqref{reg_type} with: HR estimator \eqref{HT_type}, synthetic (SYN), EBLUP and GREG estimators. The later three estimators are taken by formulas (4.2.2), (7.2.16) and (2.3.6) from \cite{R_2003}, and we denote them by $\hat t_{y;{\cal D}}^{SYN}$, $\hat t_{y;{\cal D}}^{EBLUP}$ and $\hat t_{y;{\cal D}}^{GREG}$, respectively.

The simulations are based on populations of two types. Let $N=500$ and $m=50$. The values of the variable $z$ in the population ($P_1$) are obtained as follows: for the elements of the domain ${\cal D}$, they are generated from the distribution ${\cal N}(4, 1)$, and, for the elements in ${\cal U}\backslash{\cal D}$, from ${\cal N}(6, 1.25)$. To get the values of $z$ for the population of the type ($P_2$), we generate the numbers from the exponential distribution ${\cal E}(1)$ for the whole population, and multiply them by $2$ for the elements outside the domain ${\cal D}$. Next, for each population, in order to generate and fix values of the variable $x$, we use model \eqref{x_reg_z} of the form $X_i=z_i+\delta_i$ where the independent errors $\delta_1, \ldots, \delta_N$ are distributed by ${\cal N}(0, \tau^2)$. We are interested in different correlations between $x$ and $z$, therefore, firstly, we set the variances $\tau^2$ which return the correlations $\rho_{xz}=0.1,0.2, \ldots, 0.9$ in average, and, secondly, for each population, we choose $9$ particular collections of the values of $x$ which realize the expected correlations with a small error. Next, for each set of the values of $x$, we generate collections of the values of $y$ similarly (choosing variances of a model), in order to ensure the correlations close to $\rho_{yx}=0.1,0.2, \ldots, 0.9$. But here, we use three different data generation mechanisms. 

Case (A). Model \eqref{y_reg_x} of the form $Y_i=2+x_i+\varepsilon_i$, with the independent errors $\varepsilon_i$ distributed by ${\cal N}(0, \sigma^2)$, is applied for all $i\in {\cal U}$.

Case (B). The model is $Y_i=2+\beta_{2i}x_i+\varepsilon_i$, with the independent errors $\varepsilon_i$ by ${\cal N}(0, \sigma^2)$, where $\beta_{2i}=1.25$ for $i\in {\cal D}$, and  $\beta_{2i}=1$ for $i\in {\cal U}\backslash{\cal D}$.

Case (C). The model is $Y_i=2+x_i+\varepsilon_i$, with the independent errors $\varepsilon_i$ from ${\cal N}(0, c_i\sigma^2)$, where $c_i=3$ for $i\in {\cal D}$, and  $c_i=1$ for $i\in {\cal U}\backslash{\cal D}$.

Finally, each of the types ($P_1$) and ($P_2$), mixed with one of Cases (A), (B) and (C), contains $81$ different trios of the sets of the values of $z$, $x$ and $y$ in the population. The sampling design $p(\cdot)$ is the simple random without replacement in all cases, and we take $n=75$. Then, in the domain ${\cal D}$, the expectation and the standard deviation of the sample size are $7.5$ and $\approx \hspace{-0.8mm}2.4$, respectively. To evaluate the MSEs of the estimators, we apply the Monte--Carlo (M--C) simulations by drawing independently $10^3$ samples without replacement from the populations. The probabilities $\theta_{i;{\cal D}}$, $i=1,\ldots, N$, are estimated using the M--C method too, see Appendix \ref{a:11}, with the number of replications $R=10^6$.

We assume that the survey statistician has no idea on the models in Cases (B) and (C), and, at the estimation, fits the simplest linear regressions by \eqref{y_reg_x}. This assumption is realistic because, in the domain ${\cal D}$, realizations of the sample size are too small in order to test differences between regressions. The variable $z$ is not incorporated into the estimators $\hat t_{y;{\cal D}}^{SYN}$, $\hat t_{y;{\cal D}}^{EBLUP}$ and $\hat t_{y;{\cal D}}^{GREG}$ here, see Section \ref{s:5} for an explanation. Note that the sampling design does not depend on this variable as well.

\begin{figure}[H]
\centering
\begin{minipage}{.5\textwidth}
  \centering
  \includegraphics[width=0.9\linewidth]{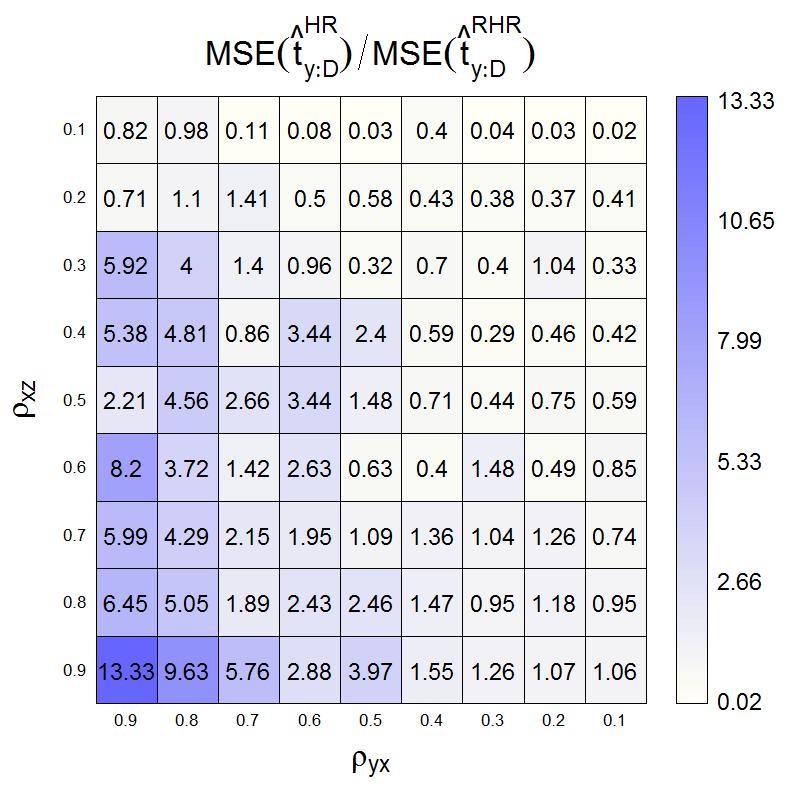}
  \captionof{figure}{($P_1$), Case (A), HR vs RHR}
  \label{fig:p2a_hr}
\end{minipage}%
\begin{minipage}{.5\textwidth}
  \centering
  \includegraphics[width=0.9\linewidth]{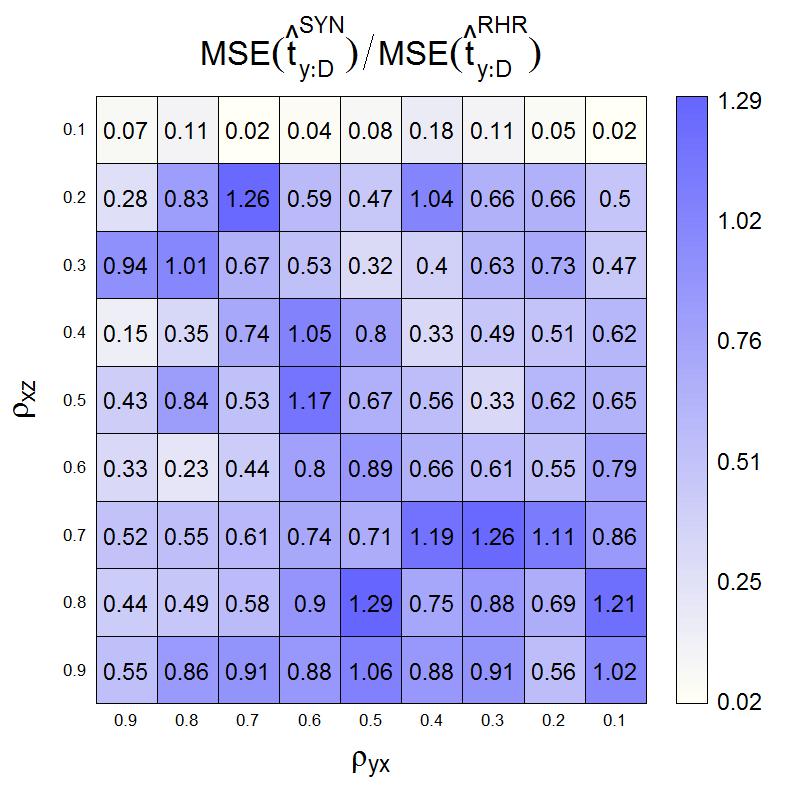}
 \captionof{figure}{($P_1$), Case (A), SYN vs RHR}
  \label{fig:p2a_syn}
\end{minipage}
\begin{minipage}{.5\textwidth}
  \centering
  \includegraphics[width=0.9\linewidth]{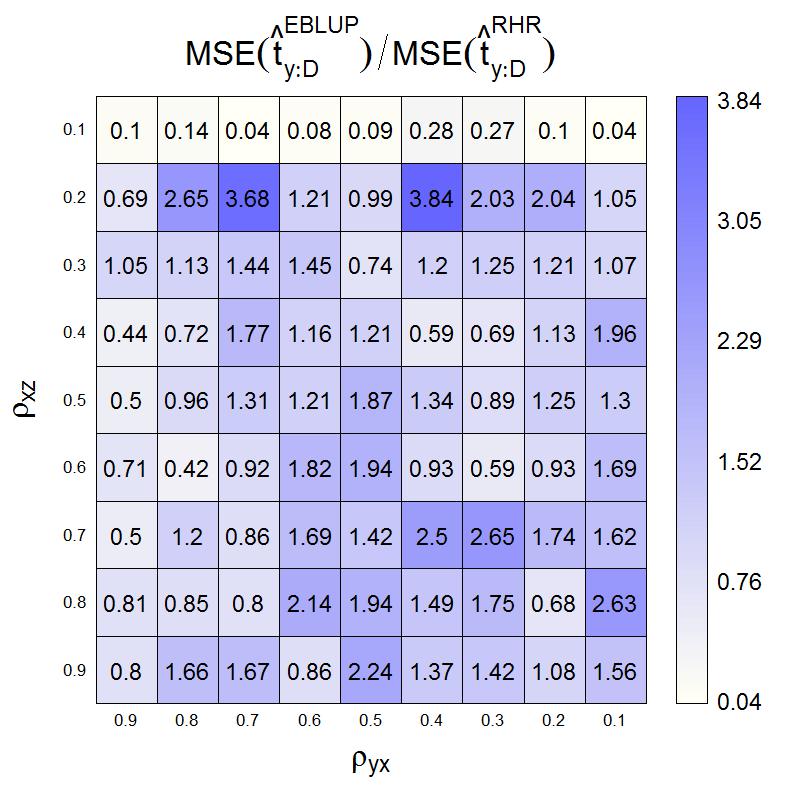}
  \captionof{figure}{($P_1$), Case (A), EBLUP vs RHR}
  \label{fig:p2a_eblup}
\end{minipage}%
\begin{minipage}{.5\textwidth}
  \centering
  \includegraphics[width=0.9\linewidth]{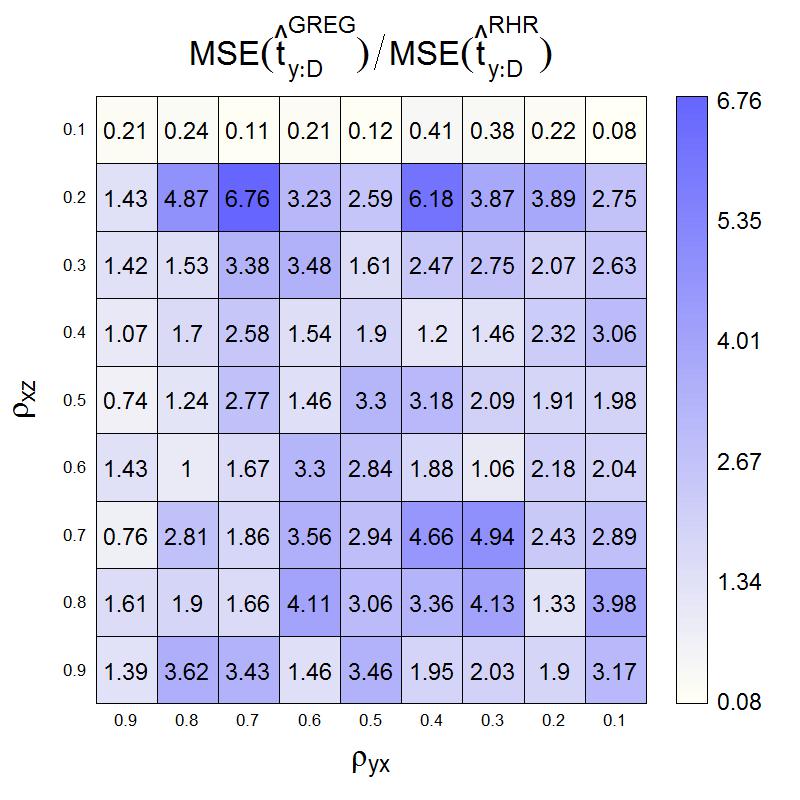}
  \captionof{figure}{($P_1$), Case (A), GREG vs RHR}
  \label{fig:p2a_greg}
\end{minipage}
\end{figure}

\begin{figure}[H]
\centering
\begin{subfigure}{.25\textwidth}
  \centering
  \includegraphics[width=0.9\linewidth]{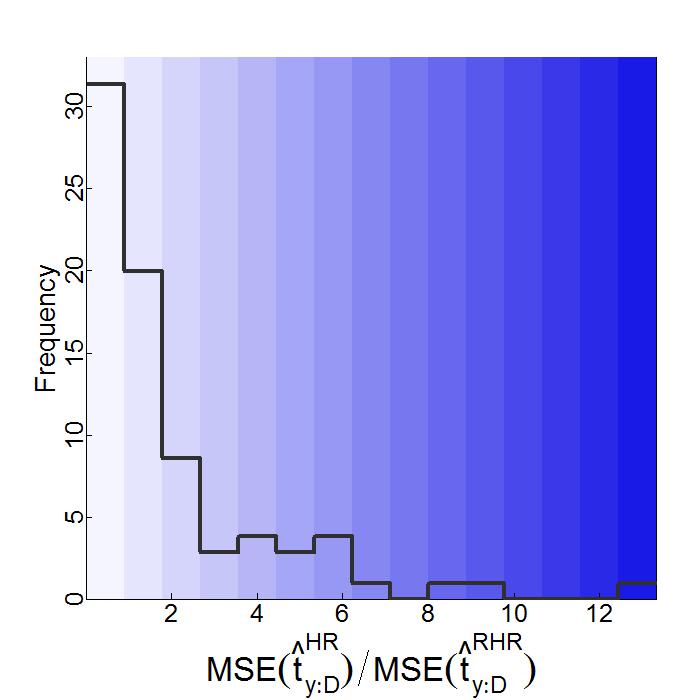}
  \label{fig:p2a_1}
\end{subfigure}%
\begin{subfigure}{.25\textwidth}
  \centering
  \includegraphics[width=0.9\linewidth]{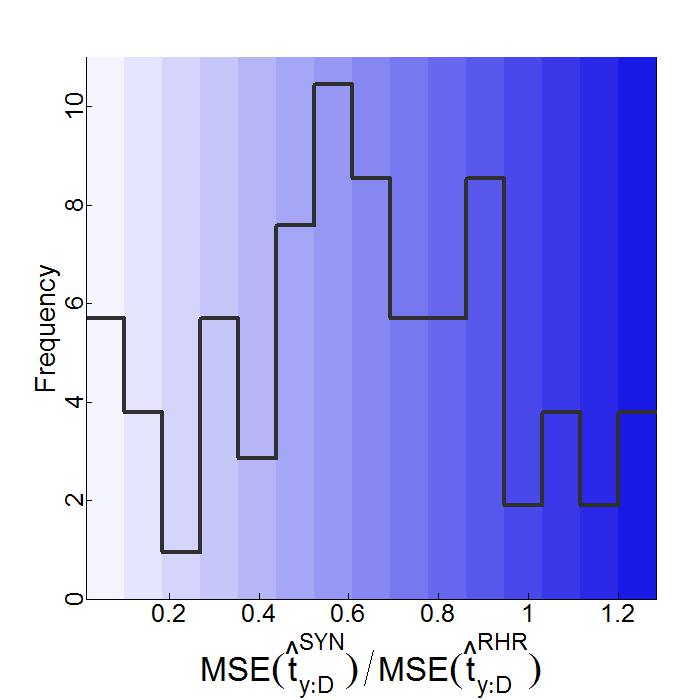}
  \label{fig:p2a_2}
\end{subfigure}%
\begin{subfigure}{.25\textwidth}
  \centering
  \includegraphics[width=0.9\linewidth]{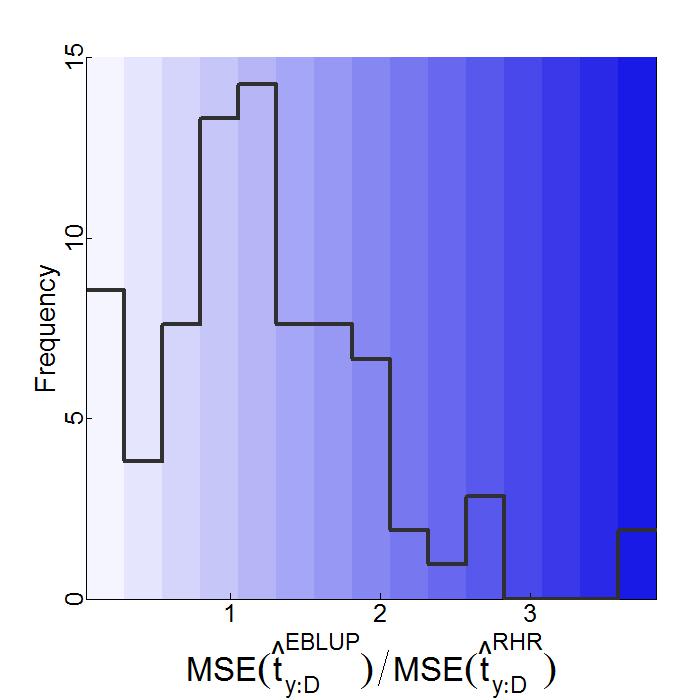}
  \label{fig:p2a_3}
\end{subfigure}%
\begin{subfigure}{.25\textwidth}
  \centering
  \includegraphics[width=0.9\linewidth]{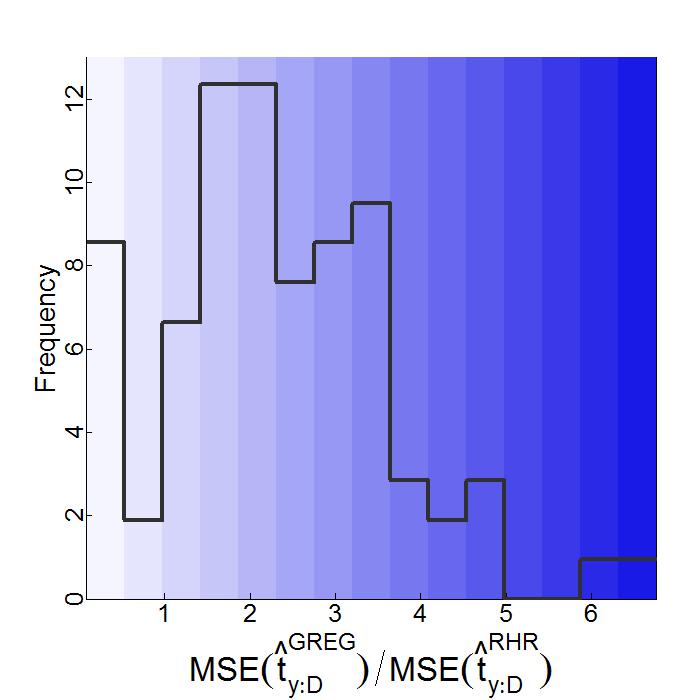}
  \label{fig:p2a_4}
\end{subfigure}
\caption{Counts for ($P_1$), Case (A), and HR, SYN, EBLUP, GREG vs RHR}
\label{fig:p2a_hist}
\end{figure}

Appendices \ref{a:21}--\ref{a:26} present simulation results on the populations ($P_1$) and ($P_2$) paired with Cases (A), (B) and (C). By the larger figures, e.g., Figures \ref{fig:p2a_hr}--\ref{fig:p2a_greg}, we present ratios of the MSEs of HR, SYN, EBLUP, GREG with the MSE of RHR, respectively. The following smaller figures, e.g., inside Figure \ref{fig:p2a_hist}, summarize the previous four. 

As we expected, the RHR estimator improves HR under stronger correlations $\rho_{yx}$ (and with stronger $\rho_{xz}$). It holds for $5$ of $6$ cases -- in Figures \ref{fig:p2a_hr}, \ref{fig:p3a_hr}, \ref{fig:p3b_hr}, \ref{fig:p2c_hr}, \ref{fig:p3c_hr} except Figure \ref{fig:p2b_hr} in Appendix \ref{a:23}. Therefore, in \ref{a:23}, we compare the HR estimator with the other ones. Let us focus further, considering figures, on more reliable correlations between the `size' variables $x$ and $z$, for instance, $\rho_{xz}\geq 0.4$. Then, in some of $5$ cases, the SYN estimator is a better predictor than RHR where the correlation $\rho_{yx}$ is strong, but it loses its power more than RHR when this correlation decreases (Figures \ref{fig:p2a_syn}, \ref{fig:p3a_syn}, \ref{fig:p3b_syn}, \ref{fig:p2c_syn}). The SYN estimator outperforms RHR in Figure \ref{fig:p2a_syn}. However, Case (C) seems the most unsuccessful for SYN (Figures \ref{fig:p2c_syn}, \ref{fig:p3c_syn}), as well as in the cases with the type ($P_2$) of the populations. In the separate case, where we compare the HR estimator with SYN, the picture (Figure \ref{fig:p2b_syn}) is different. Here one can state that, for $\rho_{xz}\geq 0.8$, the SYN estimator is better than HR, but, for $0.4\leq \rho_{xz}\leq 0.6$, the result is opposite. Similarly, but in $2$ of all $6$ cases only, the EBLUP estimator works well for strong correlations $\rho_{yx}$, see Figures \ref{fig:p2a_eblup}, \ref{fig:p2b_eblup}. In all cases of the populations ($P_2$), the RHR estimator is evidently better than EBLUP (Figures \ref{fig:p3a_eblup}, \ref{fig:p3b_eblup}, \ref{fig:p3c_eblup}). From the side of different generation models for the values of $y$, large differences in efficiencies of RHR and EBLUP are in Case (C) (Figures \ref{fig:p2c_eblup}, \ref{fig:p3c_eblup}) as well. Our estimation approach is much better than the GREG estimator except in Appendix \ref{a:23} where GREG improves the HR estimator for strong correlations $\rho_{yx}$ and $\rho_{xz}$. 

We conclude, from the simulations, that the RHR estimator competes well with the SYN, EBLUP and GREG estimators under the non-normal populations ($P_2$), and where the underlying models of Cases (B) and (C) are misspecified.

\section{Conclusions}\label{s:5}

In the simulation examples, the particular special case $X_i=z_i+\delta_i$ of model \eqref{x_reg_z}, with the independent and normally distributed errors, is taken intentionally. In this way, we show that, in order to apply estimators \eqref{HT_type} and \eqref{reg_type}, it is not necessarily to have two auxiliary variables at least. In particular, we can generate the population values of $z$ from the available values of the variable $x$ (using the corresponding regression). Then a choice of the variance $\tau^2$ of the generating model is an optimization problem. While it is unsolved, the present simulation study could suggest to set weaker $\rho_{xz}$ for estimate of $\rho_{yx}$ indicating a strong correlation, and vice versa. More specifically, with this rule, the correlation $\rho_{xz}$ could vary between $0.4$ and $0.9$. The size of the parameter $\rho_{xz}$ (or $\tau^2$) can be interpreted as a noise level which serves in a smoothing of the data for the estimation domain. If two auxiliary variables $x$ and $z$ are at the disposal, then the general methodology presented seems more natural and is much free for specific interpretations. If there are more than two auxiliary variables, generalizations of the HR method are possible.  

By the definition, HR estimator \eqref{HT_type} has the important additivity (coherence) property: if the population is partitioned into non-overlapping domains, then the sum of the HR estimators over the domains coincides with the H--T estimator of the whole population. But regression type HR estimator \eqref{reg_type} is not additive.   

Many of well-known small area estimators are expressible as a linear combination (composition) of a direct estimator and a synthetic one. According to Remark \ref{rem:1}, estimator \eqref{HT_type} is the composition too, but it is not linear. We formulate the hypothesis 
\begin{equation*}
\hat t_{y;{\cal D}}^{HR}\approx \rho_{xz}^2 \hat t_{y;{\cal D}}^{HT}+\left(1-\rho_{xz}^2\right) \hat t_{y;{\cal D}}^{S},
\end{equation*}
which we applied in the construction of estimator \eqref{bias_est} of bias \eqref{HT_bias} of the HR estimator. 

As the simulation study indicates, estimators \eqref{HT_type} and \eqref{reg_type} are comparatively robust against model \eqref{y_reg_x} misspecifications. The additivity feature implies also that the HR estimator can be efficient for domains of any size in the population.

\appendix

\section{Appendix}\label{a:1}

\subsection{Monte--Carlo estimates of $\bm{p_{ij}}$}\label{a:11}

Depending on the specification of model \eqref{x_reg_z}, firstly, we use the data $(x_i, z_i)$, $i=1,\ldots, N$ to get estimates (finite population characteristics)  $\hat\alpha_1$, $\hat\alpha_2$ and $\hat\tau_1,\ldots, \hat\tau_N$ of the superpopulation model parameters, and also to obtain estimates $\widehat F_1(\cdot), \ldots, \widehat F_N(\cdot)$. Secondly, we apply the estimated model $\tilde\eta$:
\begin{equation*}
\widetilde X_i=\hat\alpha_1+\hat\alpha_2 z_i+\tilde\delta_i, \qquad i=1,\ldots, N,
\end{equation*} 
where $\widetilde X_i$ has the distribution function $\widehat F_i(\cdot)$ with $\Ete \widetilde X_i=\hat\alpha_1+\hat\alpha_2 z_i$ and $\Varte \widetilde X_i=\hat\tau_i^2$, to generate independently the collections $(x_{1r}, \ldots, x_{Nr})$, $r=1,\ldots, R$, where $R$ is a number of M--C iterations. Finally, for each $r$, we order the numbers: $x_{(1)r}\leq \cdots \leq x_{(N)r}$, and take
\begin{equation}\label{est_p_ij}
\hat p_{ij}=\frac{1}{R} \sum_{r=1}^R \mathbb{I}\left\{ x_{ir}=x_{(j)r}\right\}, \qquad i,j=1,\ldots, N.
\end{equation}
The relative frequences \eqref{est_p_ij} are consistent estimates of probabilities \eqref{p_ij}, as $N\to \infty$, if model \eqref{x_reg_z} assumptions are sufficient to have the consistency of its parameters estimators, and  because of the law of large numbers as $R=R(N)\to \infty$.

\subsection{Approximations to $\bm{\theta_{i;{\cal D}}}$}

We propose approximations of the following form:
\begin{equation*}\label{teta_approx}
\hat\theta_{i;{\cal D}}=\hat\theta_{i;{\cal D}}(a_0)=cm\sum_{j\in {\cal D}} \mathbb{I}\left\{ |z_j-z_i|\leq a_0 \hat\alpha_2^{-1}(\hat\tau_j^2+\hat\tau_i^2)^{1/2}\right\}, \qquad i=1,\ldots, N,
\end{equation*}
where $c$ is the normalizing constant such that $\sum_{i\in {\cal U}}\hat\theta_{i;{\cal D}}=m$, and $a_0>0$ is a chosen number. Here $\hat\alpha_2>0$ and $\hat\tau_1,\ldots, \hat\tau_N$ are estimates of model \eqref{x_reg_z} parameters $\alpha_2$ and $\tau_1,\ldots, \tau_N$ obtained from the data $(x_i, z_i)$, $i=1,\ldots, N$. An optimal $a_0$ can be chosen iteratively: starting from $a_0=0.01$ and continuing with the step of a similar order, stop iterations when the sign changes less than $4$ times in the sequence $\hat\theta_{i+1;{\cal D}}(a_0)-\hat\theta_{i;{\cal D}}(a_0)$, $i=1,\ldots, N-1$.

\section{Appendix}\label{a:2}


\subsection{Case (A). HR, SYN, EBLUP, GREG vs RHR under ($P_1$)}\label{a:21}

\subsection{Case (A). HR, SYN, EBLUP, GREG vs RHR under ($P_2$)}\label{a:22}

\begin{figure}[H]
\centering
\begin{minipage}{.5\textwidth}
  \centering
  \includegraphics[width=0.9\linewidth]{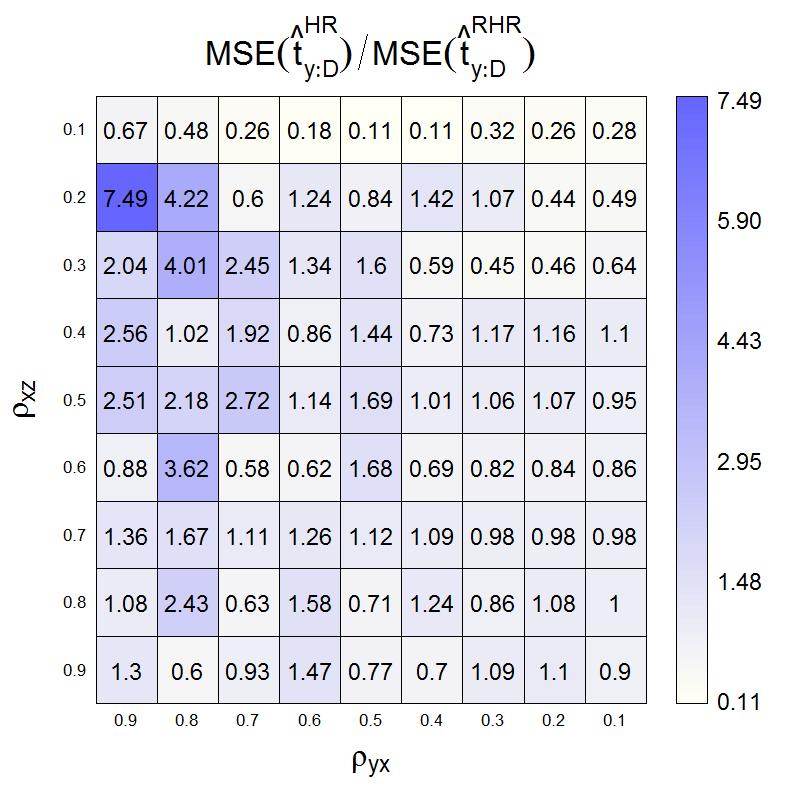}
  \captionof{figure}{($P_2$), Case (A), HR vs RHR}
  \label{fig:p3a_hr}
\end{minipage}%
\begin{minipage}{.5\textwidth}
  \centering
  \includegraphics[width=0.9\linewidth]{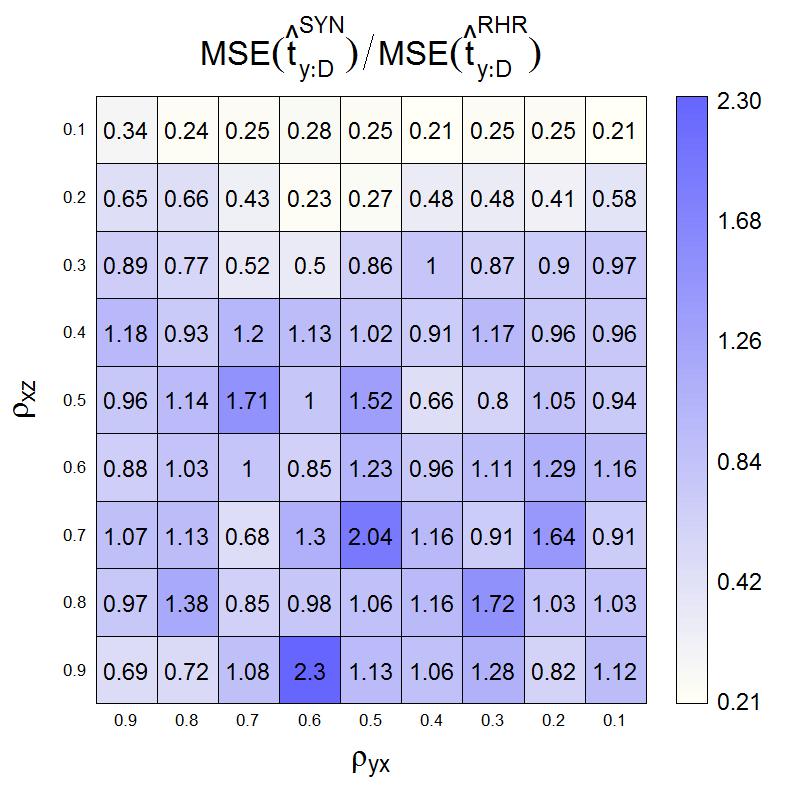}
 \captionof{figure}{($P_2$), Case (A), SYN vs RHR}
  \label{fig:p3a_syn}
\end{minipage}
\begin{minipage}{.5\textwidth}
  \centering
  \includegraphics[width=0.9\linewidth]{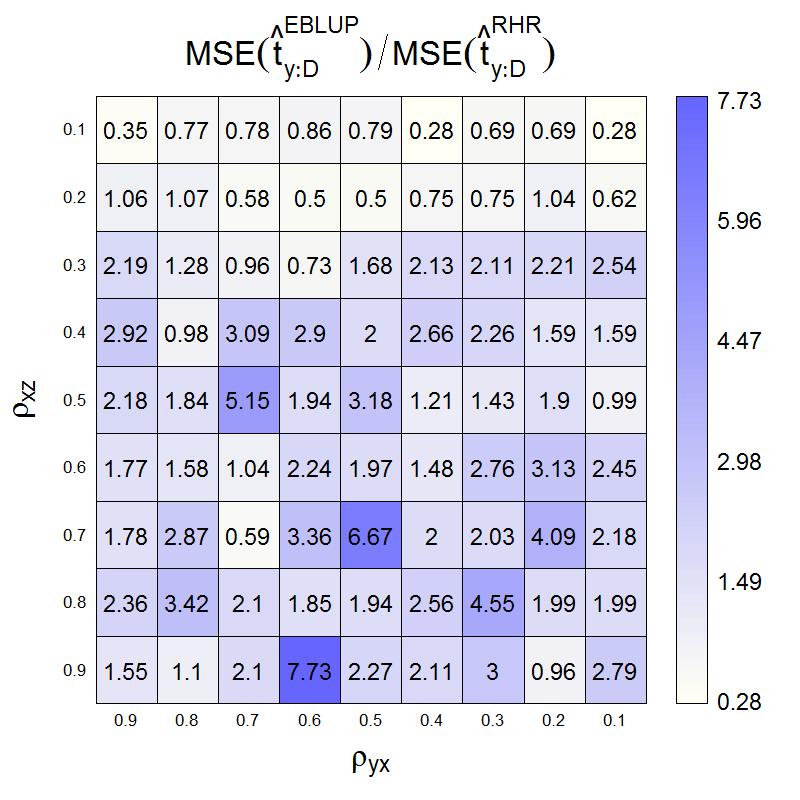}
  \captionof{figure}{($P_2$), Case (A), EBLUP vs RHR}
  \label{fig:p3a_eblup}
\end{minipage}%
\begin{minipage}{.5\textwidth}
  \centering
  \includegraphics[width=0.9\linewidth]{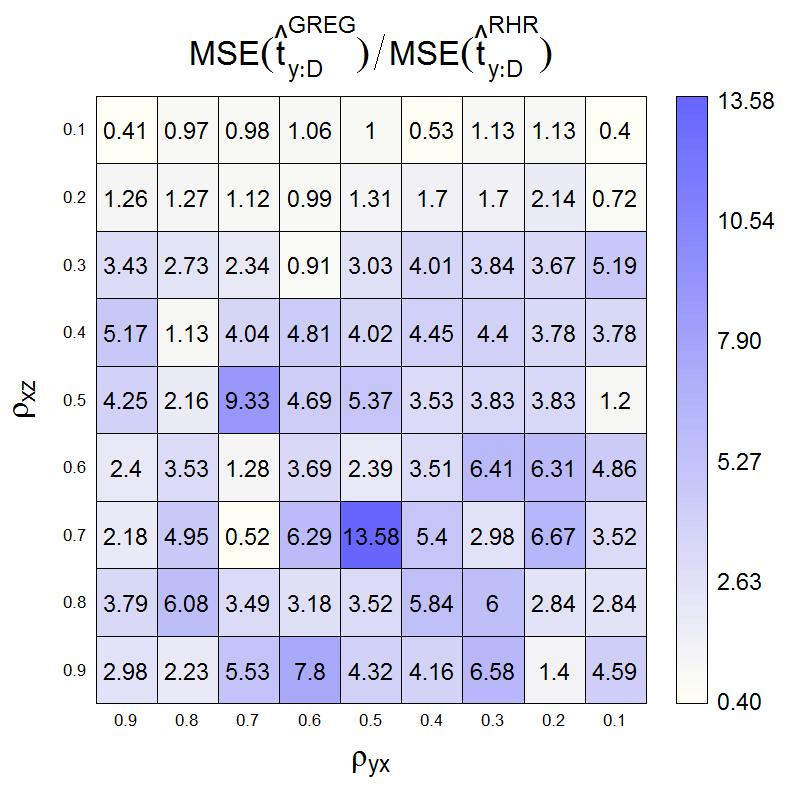}
  \captionof{figure}{($P_2$), Case (A), GREG vs RHR}
  \label{fig:p3a_greg}
\end{minipage}
\end{figure}

\begin{figure}[H]
\centering
\begin{subfigure}{.25\textwidth}
  \centering
  \includegraphics[width=0.9\linewidth]{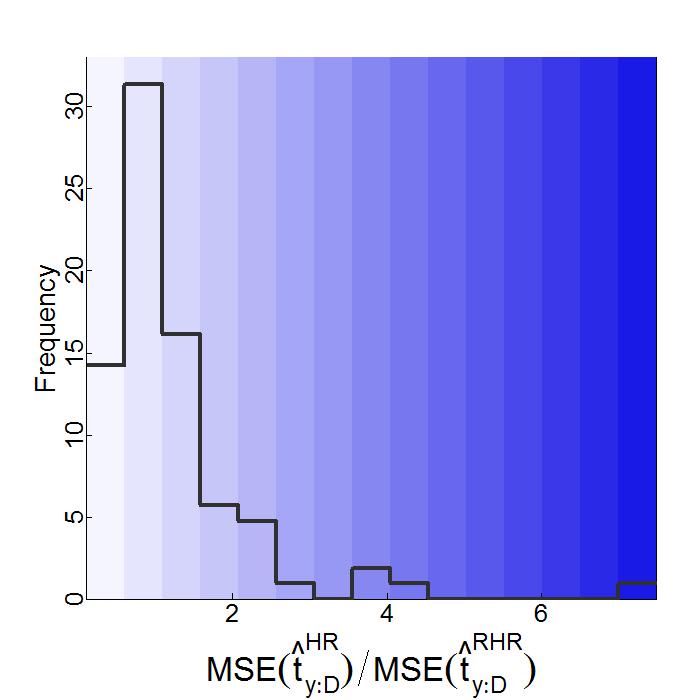}
  \label{fig:p3a_1}
\end{subfigure}%
\begin{subfigure}{.25\textwidth}
  \centering
  \includegraphics[width=0.9\linewidth]{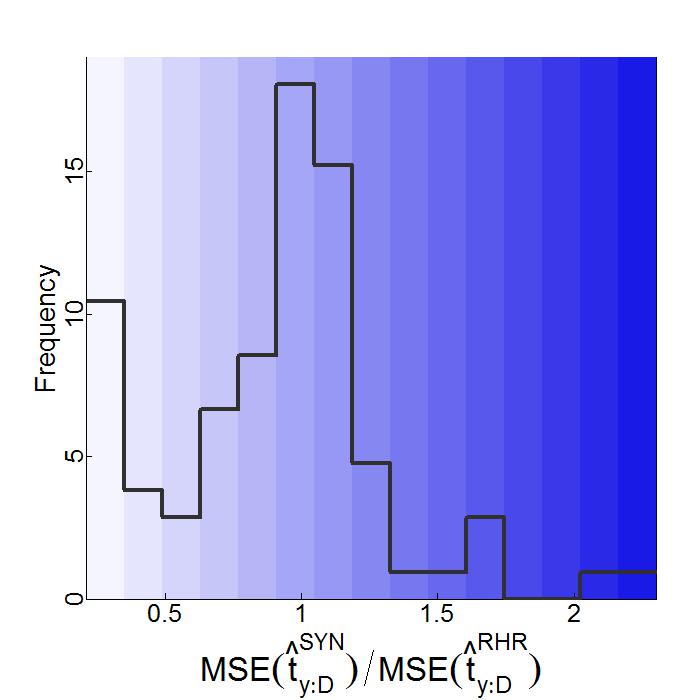}
  \label{fig:p3a_2}
\end{subfigure}%
\begin{subfigure}{.25\textwidth}
  \centering
  \includegraphics[width=0.9\linewidth]{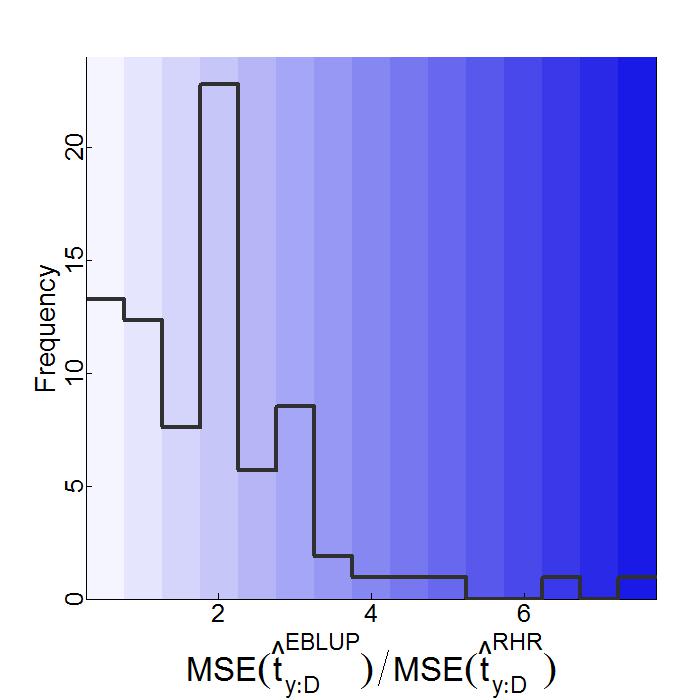}
  \label{fig:p3a_3}
\end{subfigure}%
\begin{subfigure}{.25\textwidth}
  \centering
  \includegraphics[width=0.9\linewidth]{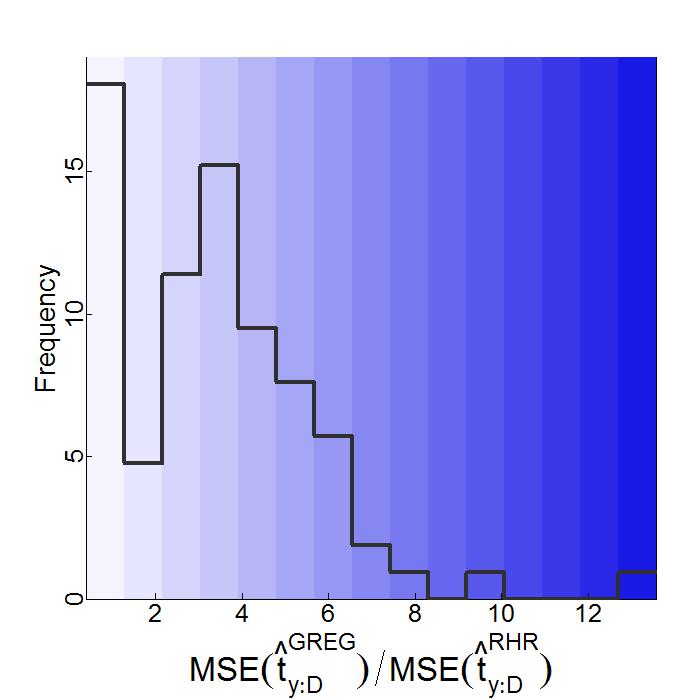}
  \label{fig:p3a_4}
\end{subfigure}
\caption{Counts for ($P_2$), Case (A), and HR, SYN, EBLUP, GREG vs RHR}
\label{fig:p3a_hist}
\end{figure}


\subsection{Case (B). RHR, SYN, EBLUP, GREG vs HR under ($P_1$)} \label{a:23}


\begin{figure}[H]
\centering
\begin{minipage}{.5\textwidth}
  \centering
  \includegraphics[width=0.9\linewidth]{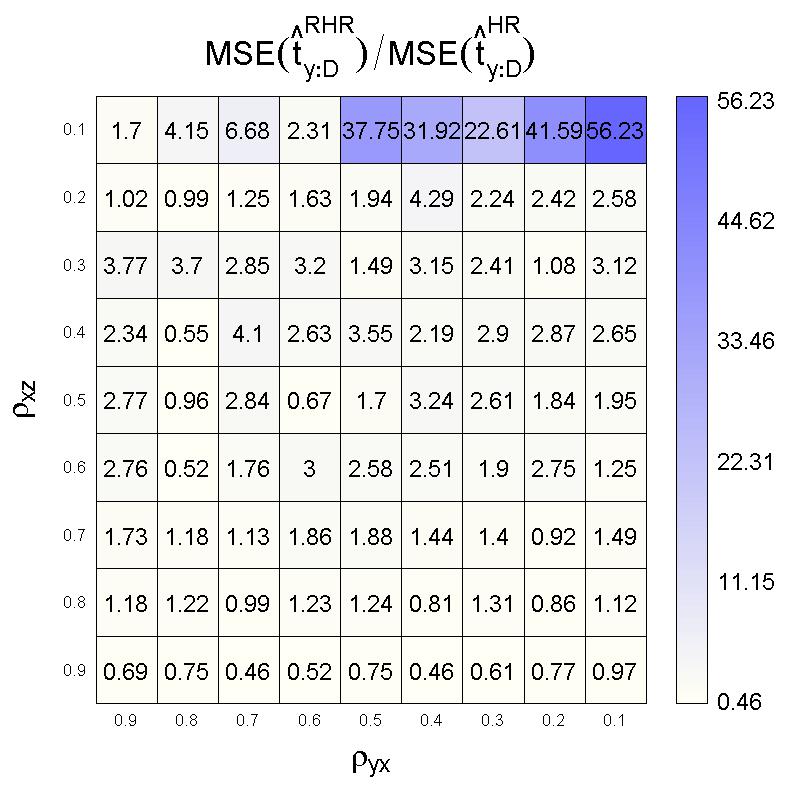}
  \captionof{figure}{($P_1$), Case (B), RHR vs HR}
  \label{fig:p2b_hr}
\end{minipage}%
\begin{minipage}{.5\textwidth}
  \centering
  \includegraphics[width=0.9\linewidth]{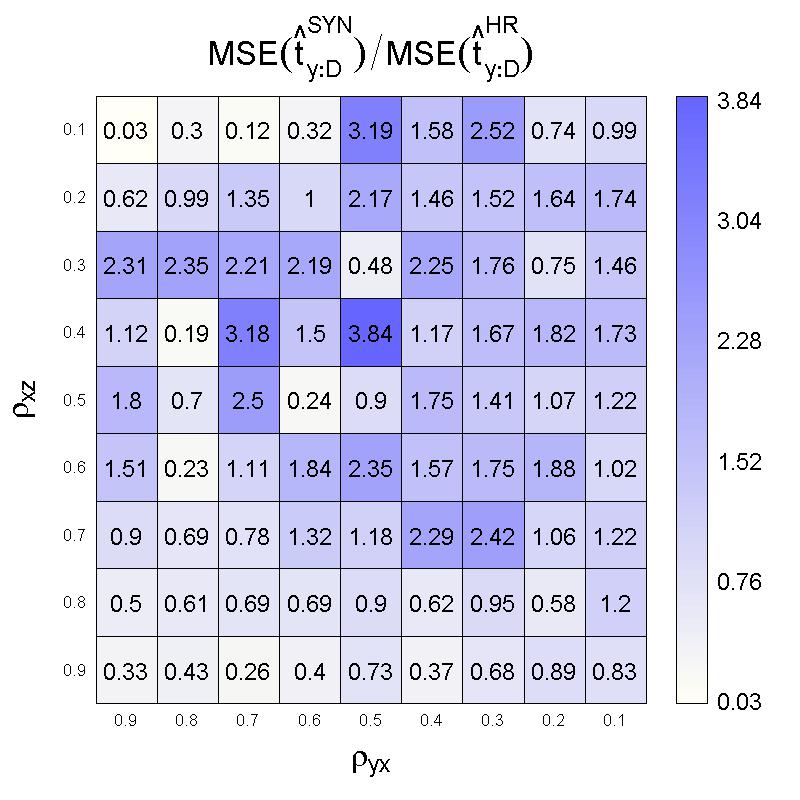}
 \captionof{figure}{($P_1$), Case (B), SYN vs HR}
  \label{fig:p2b_syn}
\end{minipage}
\begin{minipage}{.5\textwidth}
  \centering
  \includegraphics[width=0.9\linewidth]{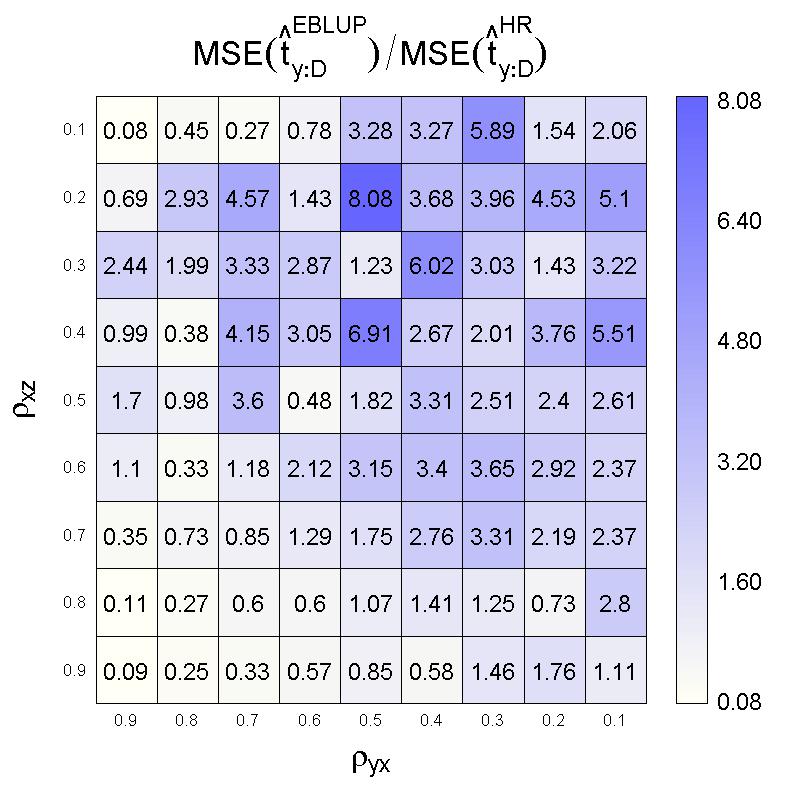}
  \captionof{figure}{($P_1$), Case (B), EBLUP vs HR}
  \label{fig:p2b_eblup}
\end{minipage}%
\begin{minipage}{.5\textwidth}
  \centering
  \includegraphics[width=0.9\linewidth]{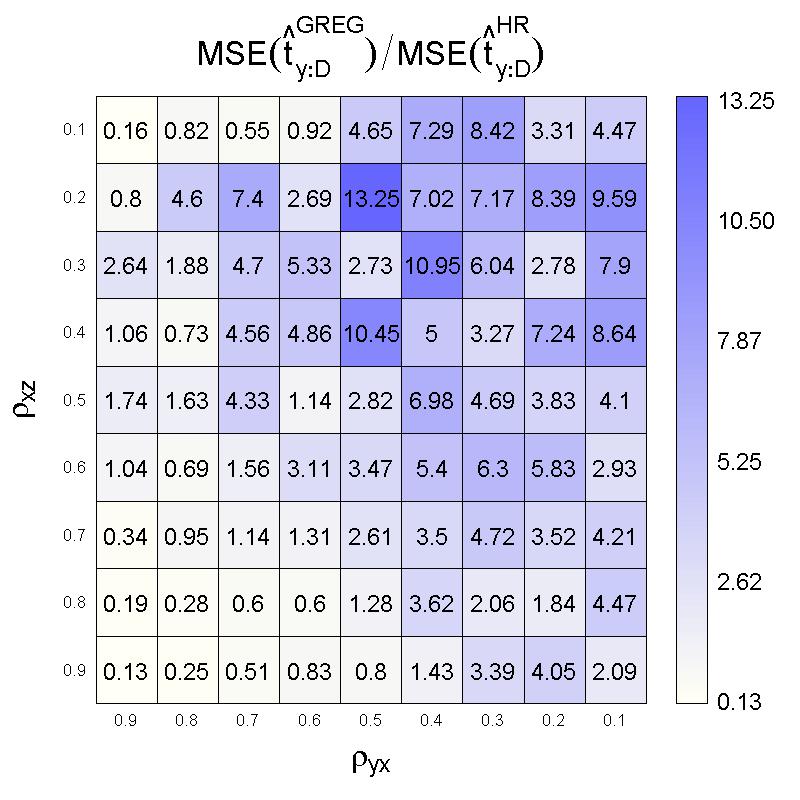}
  \captionof{figure}{($P_1$), Case (B), GREG vs HR}
  \label{fig:p2b_greg}
\end{minipage}
\end{figure}

\begin{figure}[H]
\centering
\begin{subfigure}{.25\textwidth}
  \centering
  \includegraphics[width=0.9\linewidth]{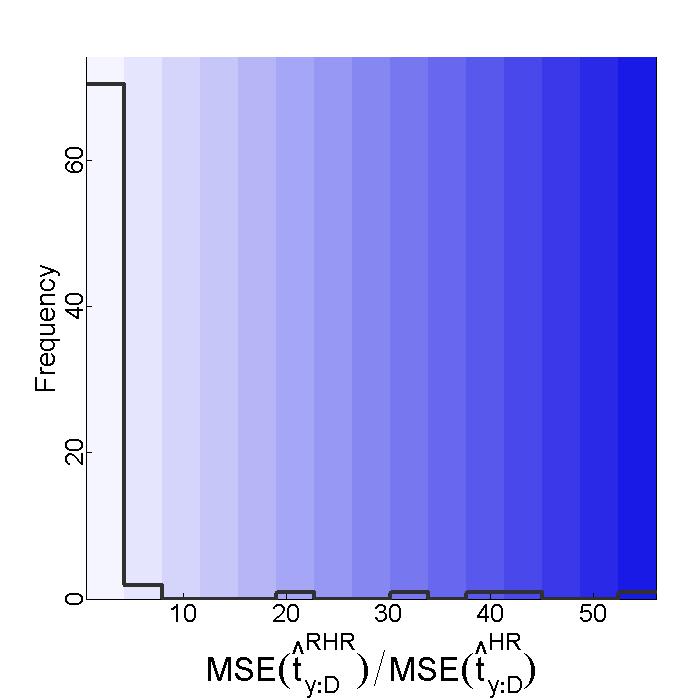}
  \label{fig:p2b_1}
\end{subfigure}%
\begin{subfigure}{.25\textwidth}
  \centering
  \includegraphics[width=0.9\linewidth]{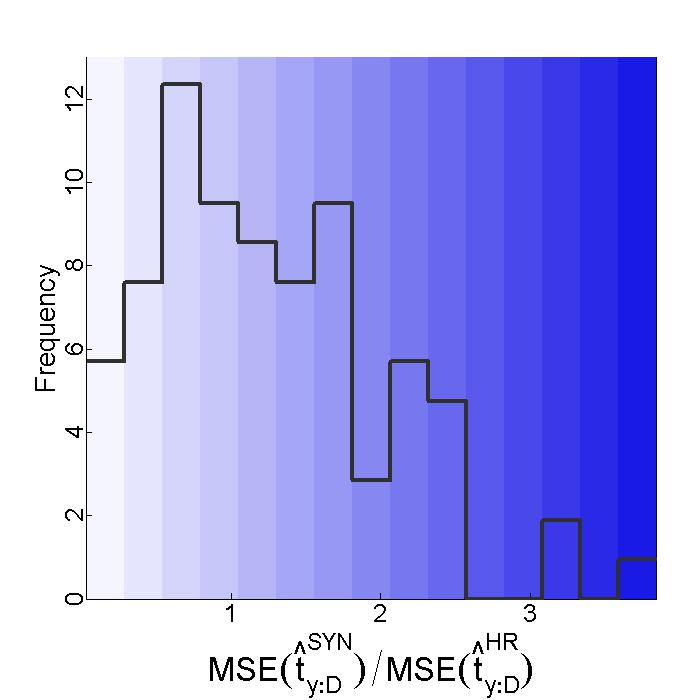}
  \label{fig:p2b_2}
\end{subfigure}%
\begin{subfigure}{.25\textwidth}
  \centering
  \includegraphics[width=0.9\linewidth]{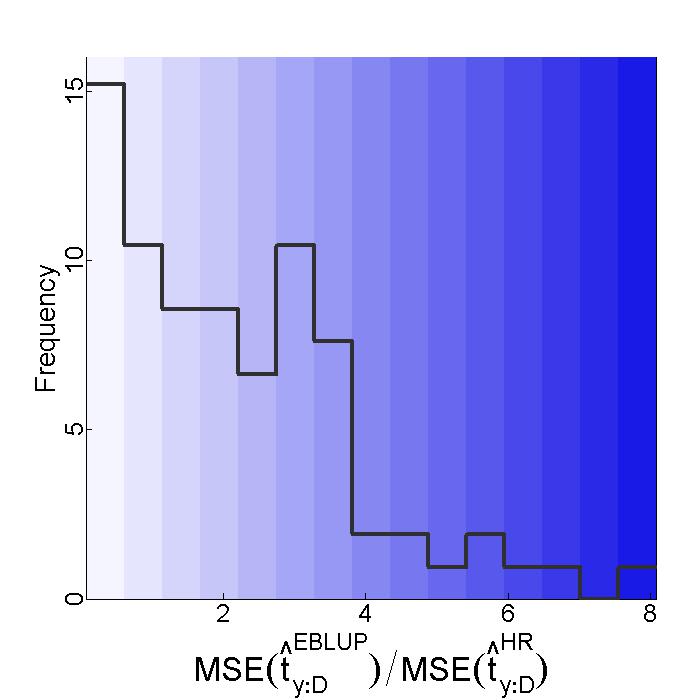}
  \label{fig:p2b_3}
\end{subfigure}%
\begin{subfigure}{.25\textwidth}
  \centering
  \includegraphics[width=0.9\linewidth]{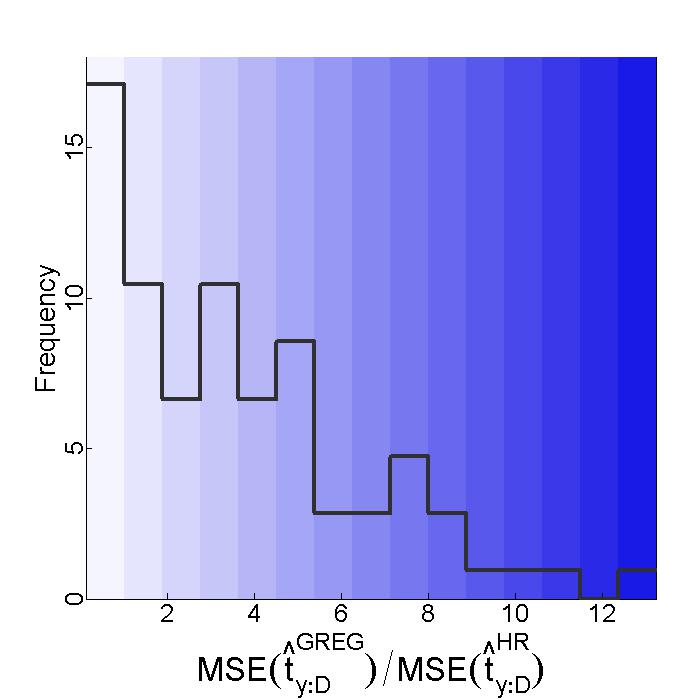}
  \label{fig:p2b_4}
\end{subfigure}
\caption{Counts for ($P_1$), Case (B), and RHR, SYN, EBLUP, GREG vs HR}
\label{fig:p2b_hist}
\end{figure}


\subsection{Case (B). HR, SYN, EBLUP, GREG vs RHR under ($P_2$)} \label{a:24}

\begin{figure}[H]
\centering
\begin{minipage}{.5\textwidth}
  \centering
  \includegraphics[width=0.9\linewidth]{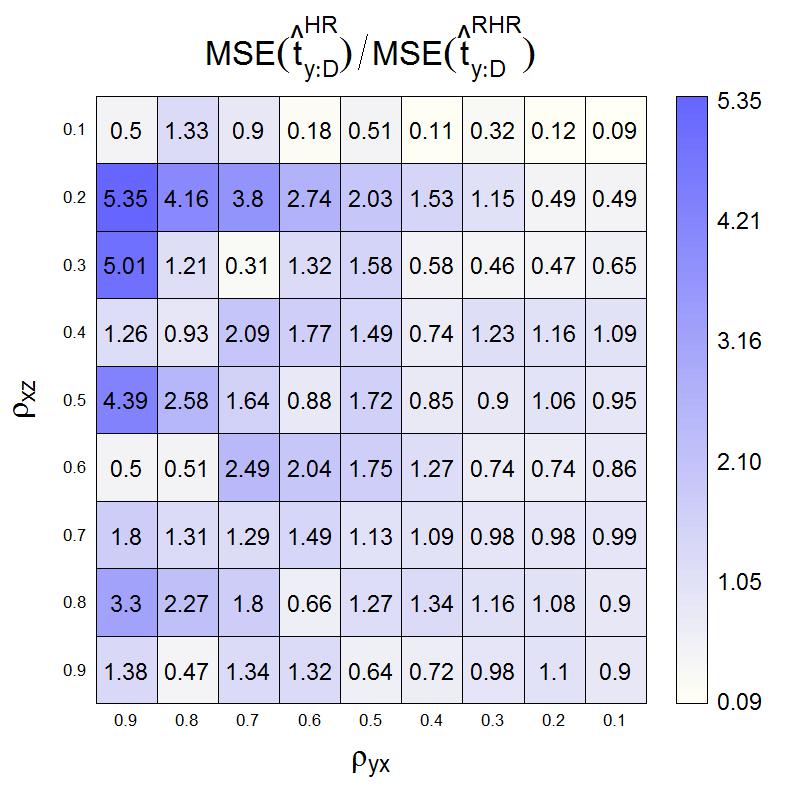}
  \captionof{figure}{($P_2$), Case (B), HR vs RHR}
  \label{fig:p3b_hr}
\end{minipage}%
\begin{minipage}{.5\textwidth}
  \centering
  \includegraphics[width=0.9\linewidth]{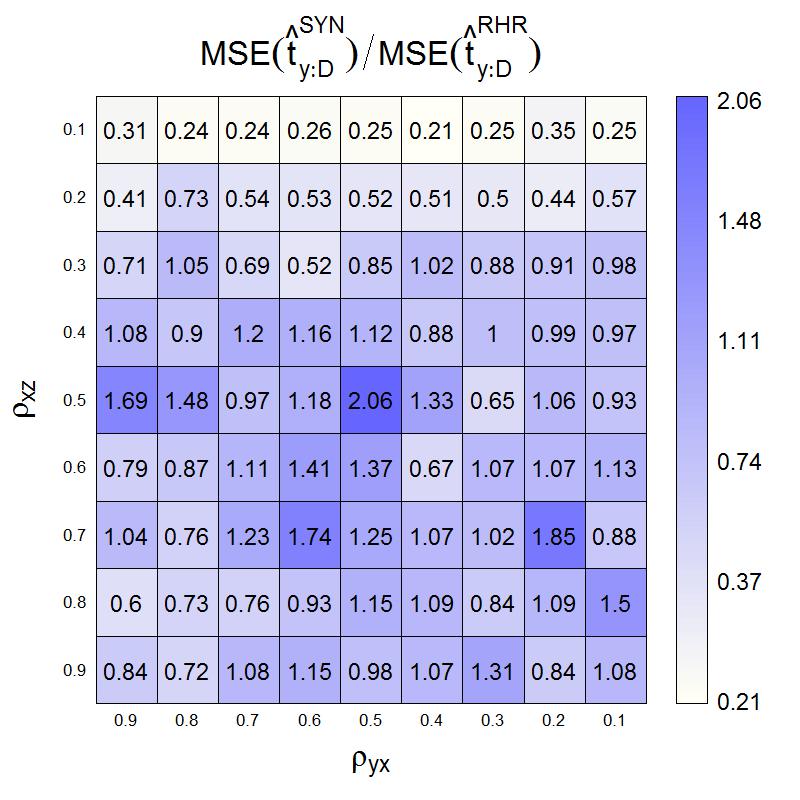}
 \captionof{figure}{($P_2$), Case (B), SYN vs RHR}
  \label{fig:p3b_syn}
\end{minipage}
\begin{minipage}{.5\textwidth}
  \centering
  \includegraphics[width=0.9\linewidth]{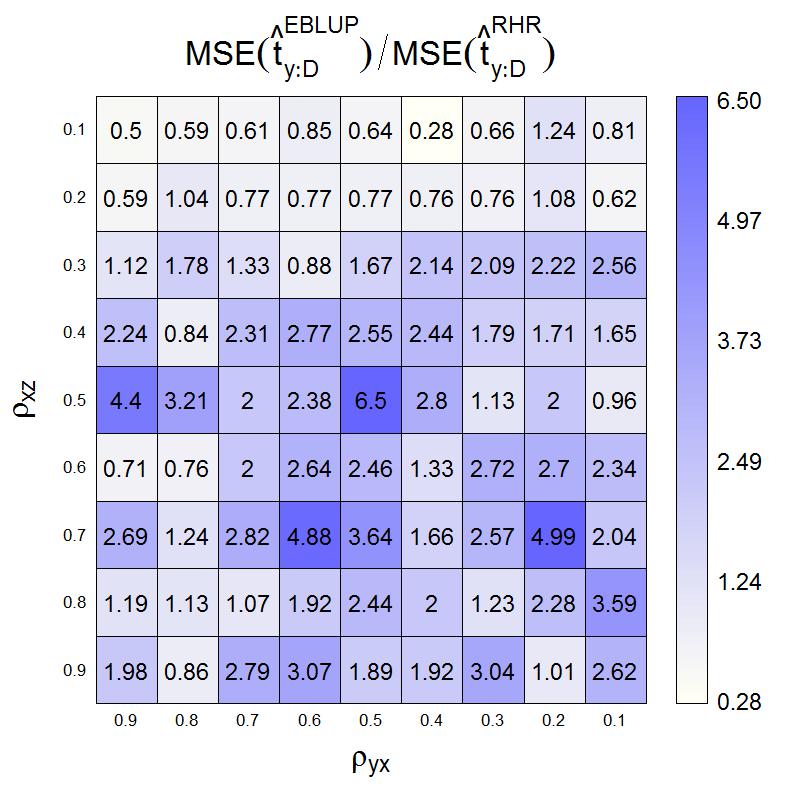}
  \captionof{figure}{($P_2$), Case (B), EBLUP vs RHR}
  \label{fig:p3b_eblup}
\end{minipage}%
\begin{minipage}{.5\textwidth}
  \centering
  \includegraphics[width=0.9\linewidth]{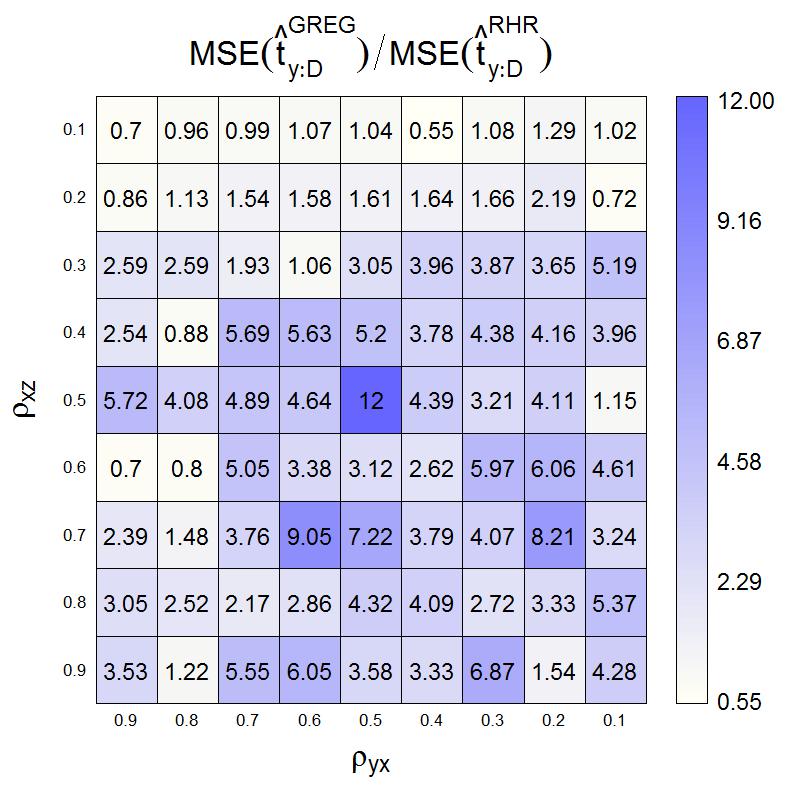}
  \captionof{figure}{($P_2$), Case (B), GREG vs RHR}
  \label{fig:p3b_greg}
\end{minipage}
\end{figure}

\begin{figure}[H]
\centering
\begin{subfigure}{.25\textwidth}
  \centering
  \includegraphics[width=0.9\linewidth]{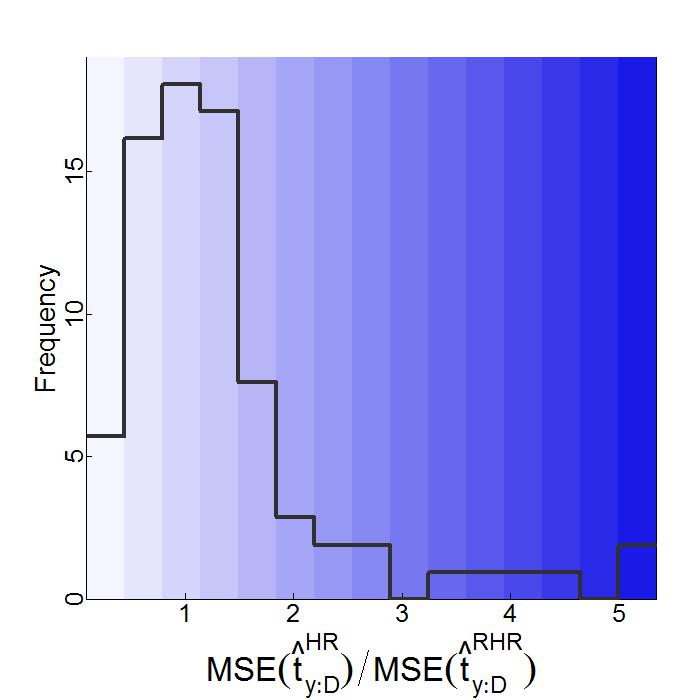}
  \label{fig:p3b_1}
\end{subfigure}%
\begin{subfigure}{.25\textwidth}
  \centering
  \includegraphics[width=0.9\linewidth]{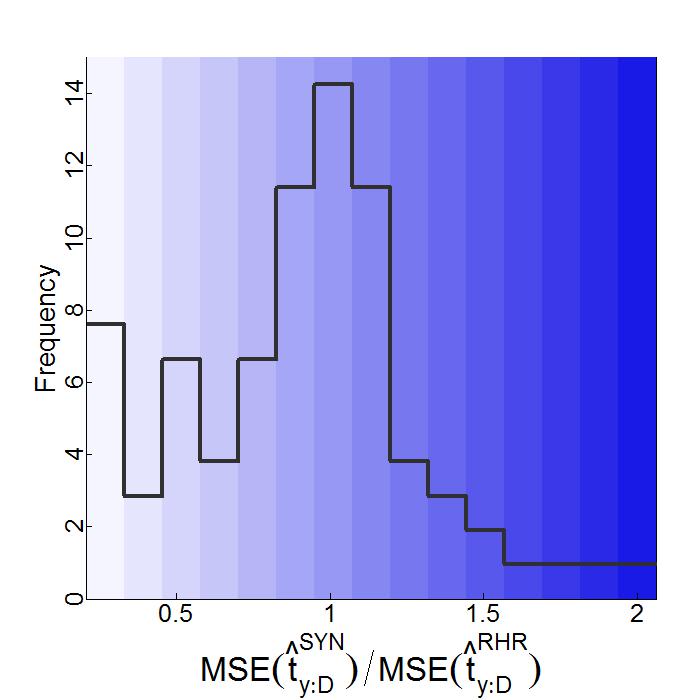}
  \label{fig:p3b_2}
\end{subfigure}%
\begin{subfigure}{.25\textwidth}
  \centering
  \includegraphics[width=0.9\linewidth]{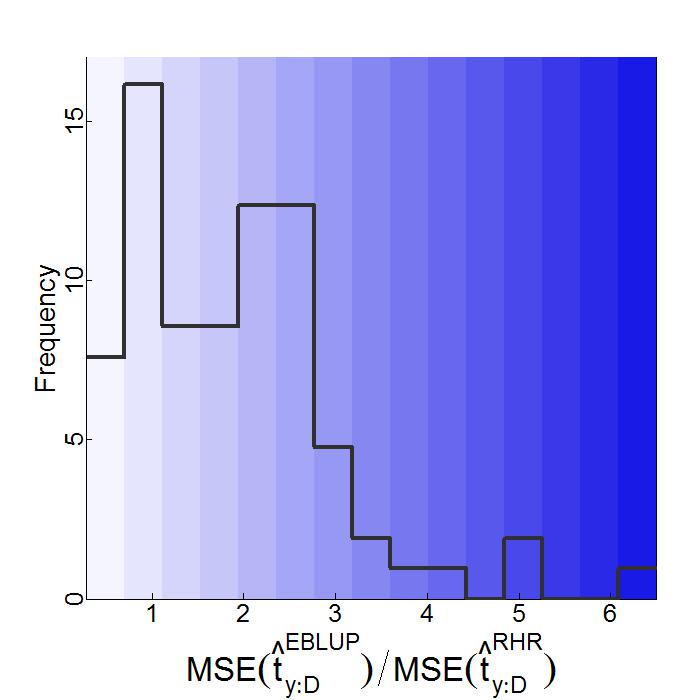}
  \label{fig:p3b_3}
\end{subfigure}%
\begin{subfigure}{.25\textwidth}
  \centering
  \includegraphics[width=0.9\linewidth]{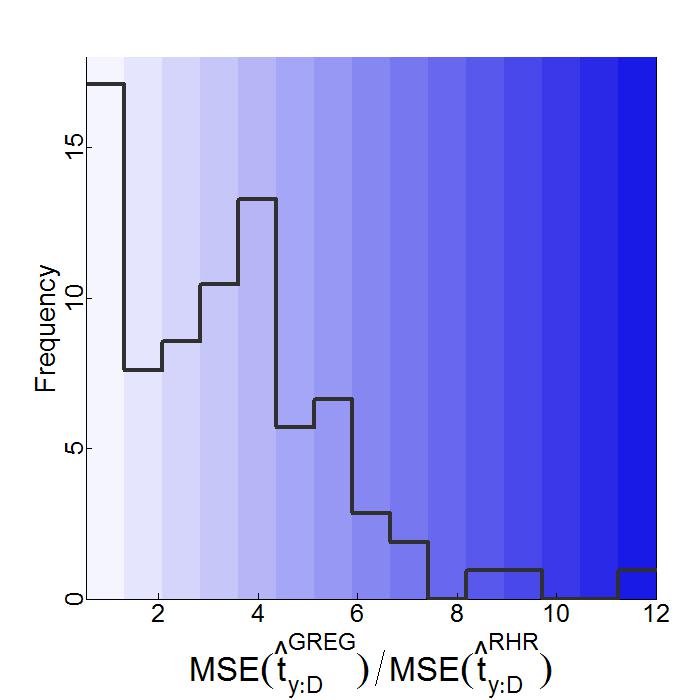}
  \label{fig:p3b_4}
\end{subfigure}
\caption{Counts for ($P_2$), Case (B), and HR, SYN, EBLUP, GREG vs RHR}
\label{fig:p3b_hist}
\end{figure}


\subsection{Case (C). HR, SYN, EBLUP, GREG vs RHR under ($P_1$)}\label{a:25}

\begin{figure}[H]
\centering
\begin{minipage}{.5\textwidth}
  \centering
  \includegraphics[width=0.9\linewidth]{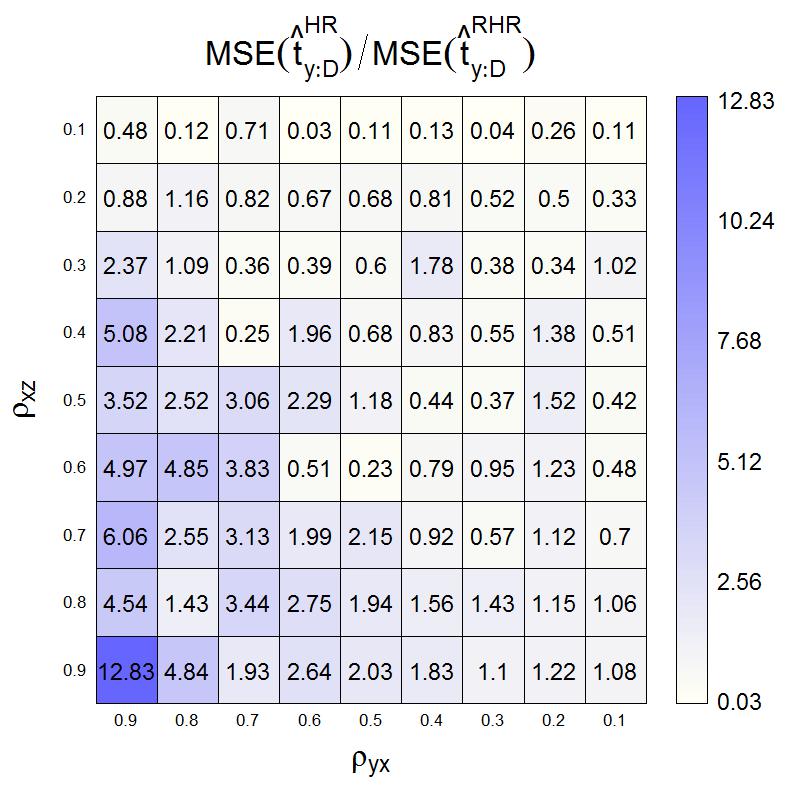}
  \captionof{figure}{($P_1$), Case (C), HR vs RHR}
  \label{fig:p2c_hr}
\end{minipage}%
\begin{minipage}{.5\textwidth}
  \centering
  \includegraphics[width=0.9\linewidth]{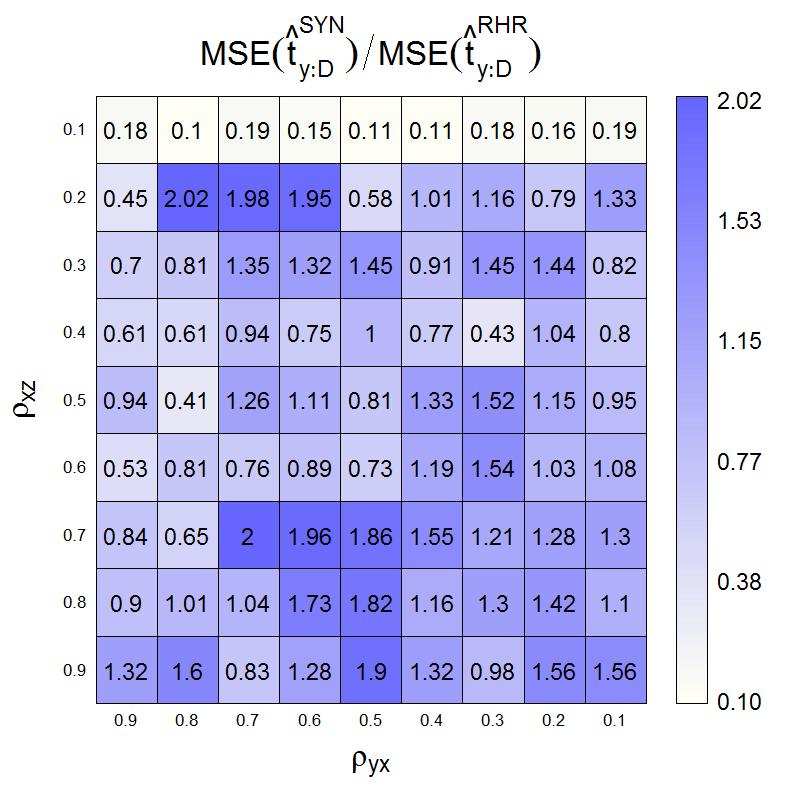}
 \captionof{figure}{($P_1$), Case (C), SYN vs RHR}
  \label{fig:p2c_syn}
\end{minipage}
\begin{minipage}{.5\textwidth}
  \centering
  \includegraphics[width=0.9\linewidth]{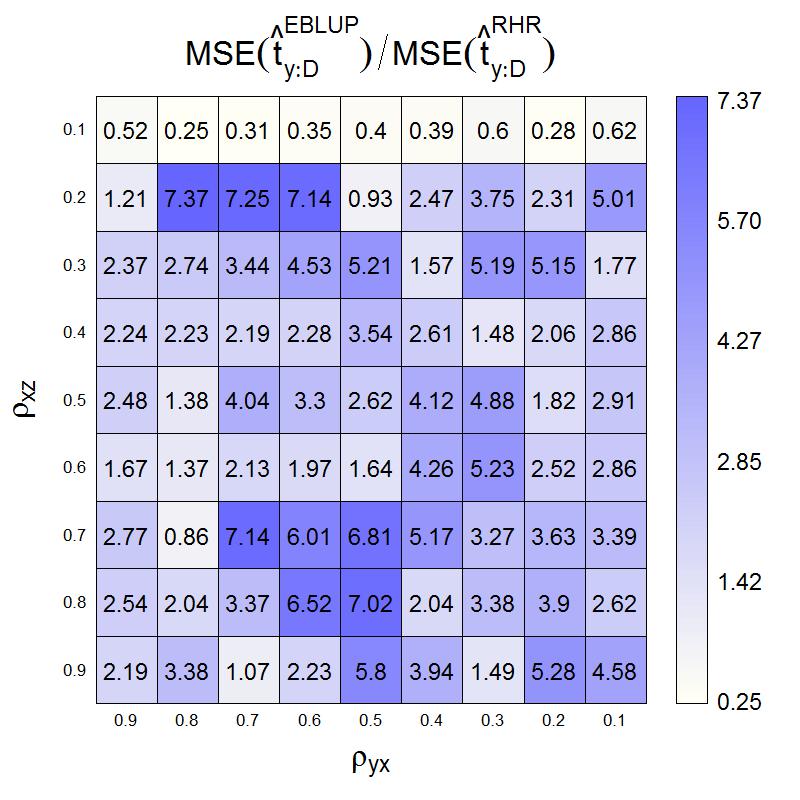}
  \captionof{figure}{($P_1$), Case (C), EBLUP vs RHR}
  \label{fig:p2c_eblup}
\end{minipage}%
\begin{minipage}{.5\textwidth}
  \centering
  \includegraphics[width=0.9\linewidth]{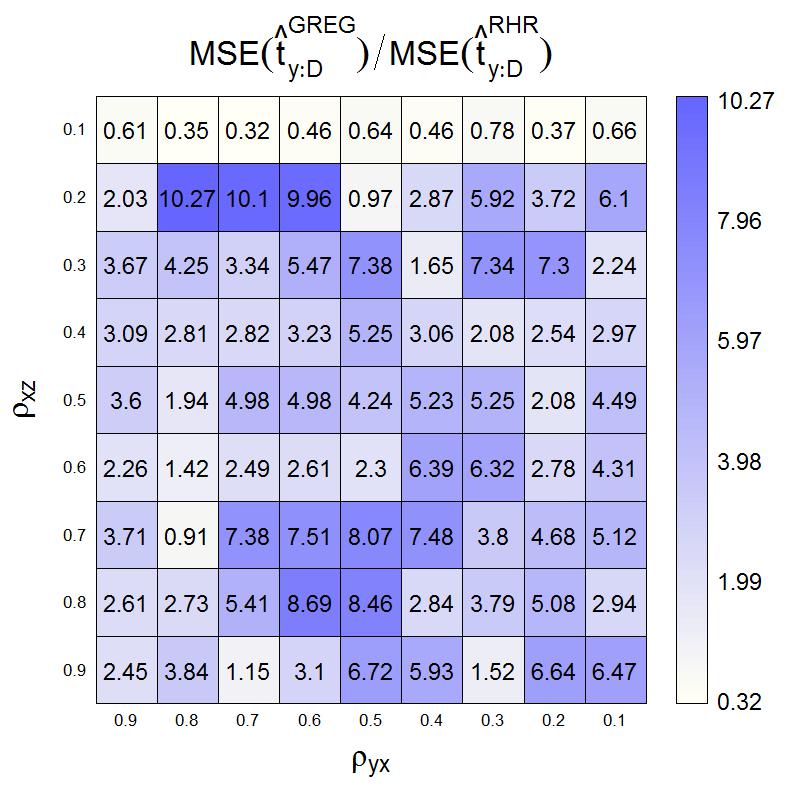}
  \captionof{figure}{($P_1$), Case (C), GREG vs RHR}
  \label{fig:p2c_greg}
\end{minipage}
\end{figure}

\begin{figure}[H]
\centering
\begin{subfigure}{.25\textwidth}
  \centering
  \includegraphics[width=0.9\linewidth]{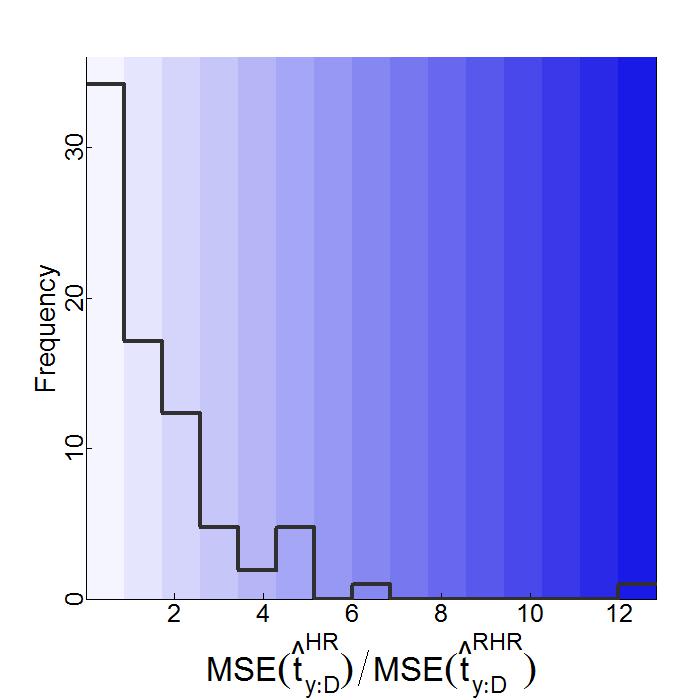}
  \label{fig:p2c_1}
\end{subfigure}%
\begin{subfigure}{.25\textwidth}
  \centering
  \includegraphics[width=0.9\linewidth]{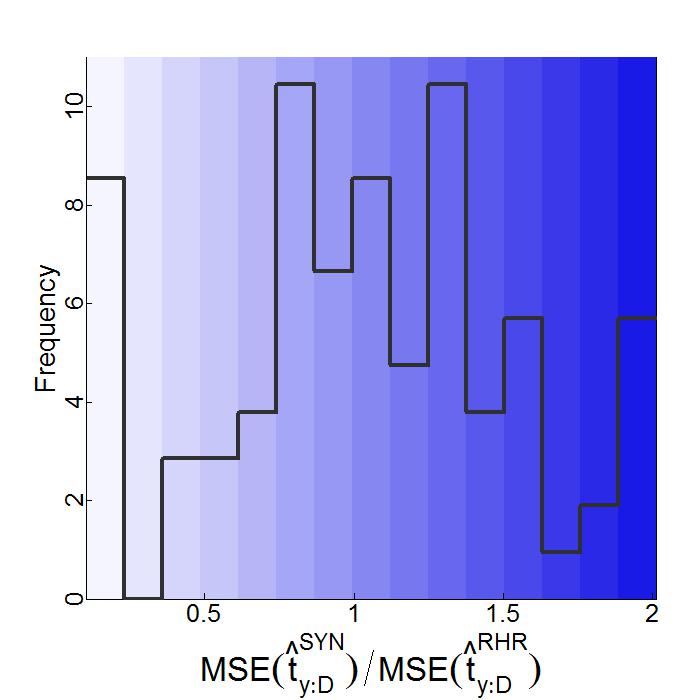}
  \label{fig:p2c_2}
\end{subfigure}%
\begin{subfigure}{.25\textwidth}
  \centering
  \includegraphics[width=0.9\linewidth]{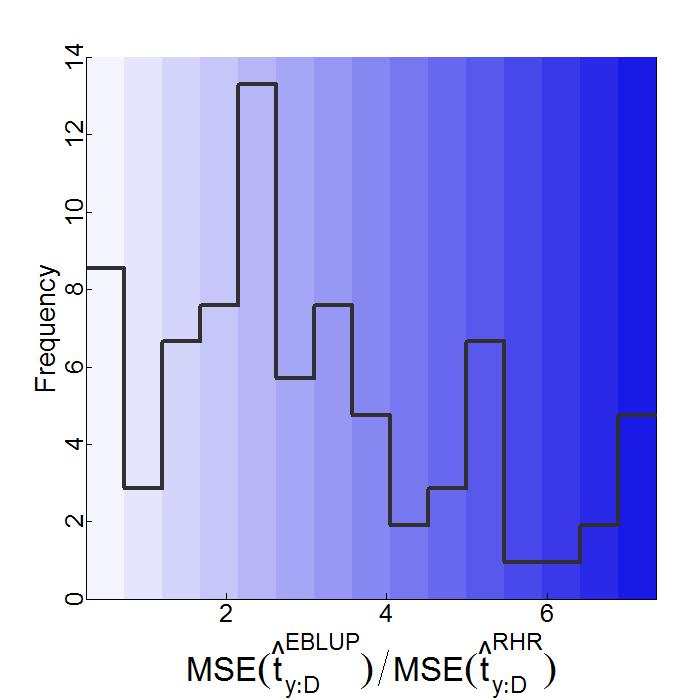}
  \label{fig:p2c_3}
\end{subfigure}%
\begin{subfigure}{.25\textwidth}
  \centering
  \includegraphics[width=0.9\linewidth]{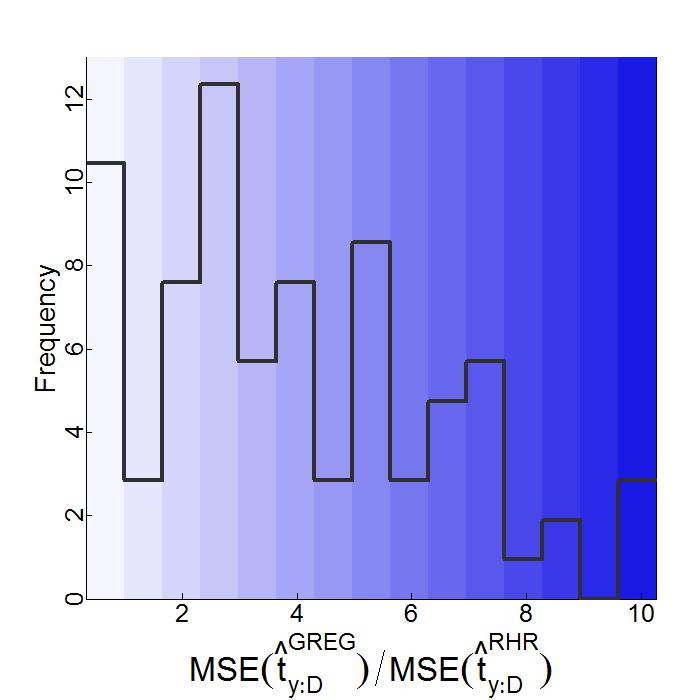}
  \label{fig:p2c_4}
\end{subfigure}
\caption{Counts for ($P_1$), Case (C), and HR, SYN, EBLUP, GREG vs RHR}
\label{fig:p2c_hist}
\end{figure}


\subsection{Case (C). HR, SYN, EBLUP, GREG vs RHR under ($P_2$)}\label{a:26}

\begin{figure}[H]
\centering
\begin{minipage}{.5\textwidth}
  \centering
  \includegraphics[width=0.9\linewidth]{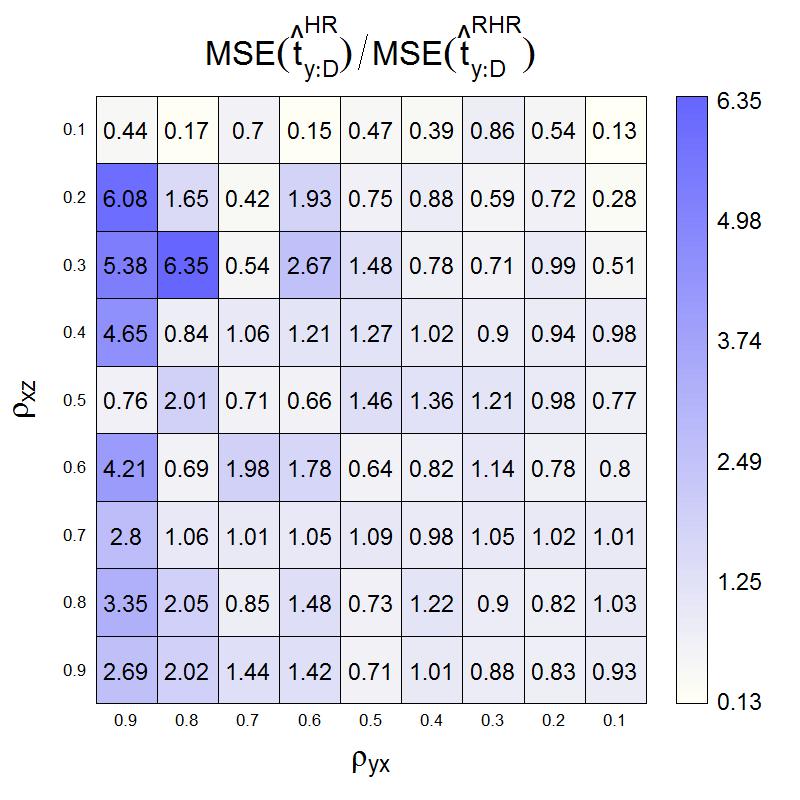}
  \captionof{figure}{($P_2$), Case (C), HR vs RHR}
  \label{fig:p3c_hr}
\end{minipage}%
\begin{minipage}{.5\textwidth}
  \centering
  \includegraphics[width=0.9\linewidth]{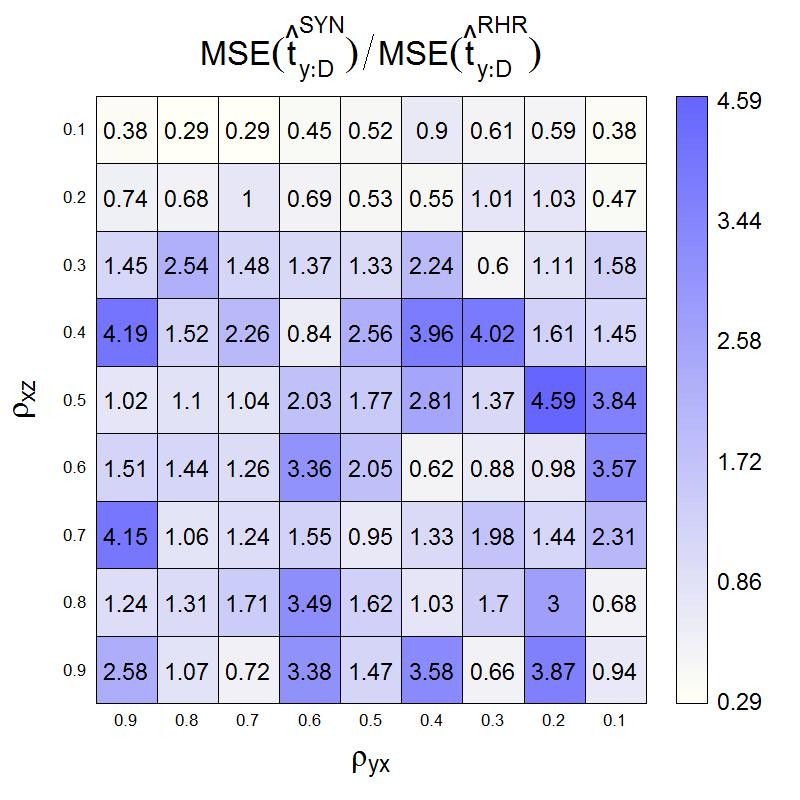}
 \captionof{figure}{($P_2$), Case (C), SYN vs RHR}
  \label{fig:p3c_syn}
\end{minipage}
\begin{minipage}{.5\textwidth}
  \centering
  \includegraphics[width=0.9\linewidth]{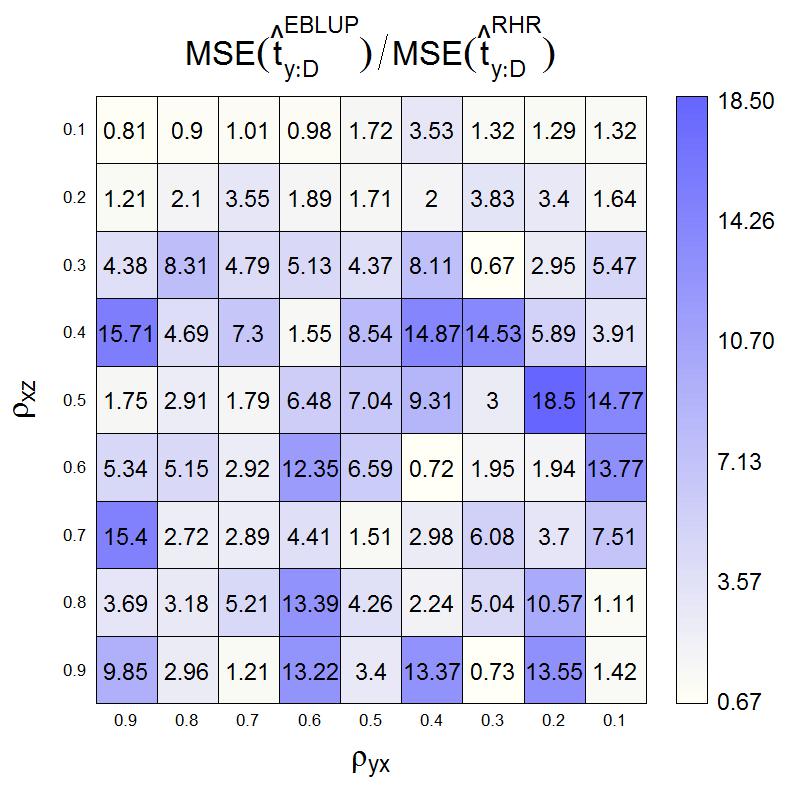}
  \captionof{figure}{($P_2$), Case (C), EBLUP vs RHR}
  \label{fig:p3c_eblup}
\end{minipage}%
\begin{minipage}{.5\textwidth}
  \centering
  \includegraphics[width=0.9\linewidth]{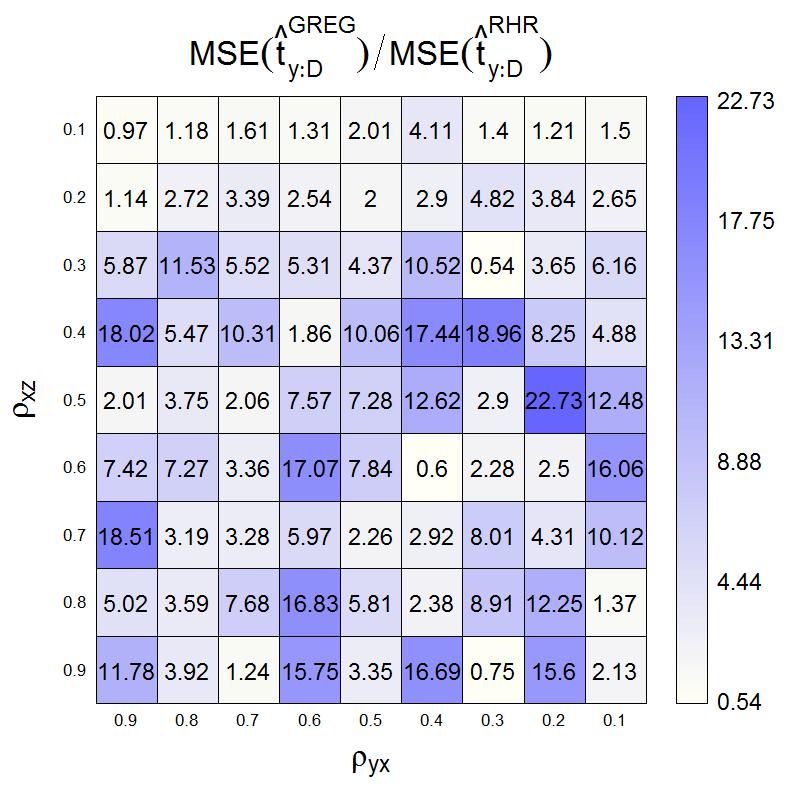}
  \captionof{figure}{($P_2$), Case (C), GREG vs RHR}
  \label{fig:p3c_greg}
\end{minipage}
\end{figure}

\begin{figure}[H]
\centering
\begin{subfigure}{.25\textwidth}
  \centering
  \includegraphics[width=0.9\linewidth]{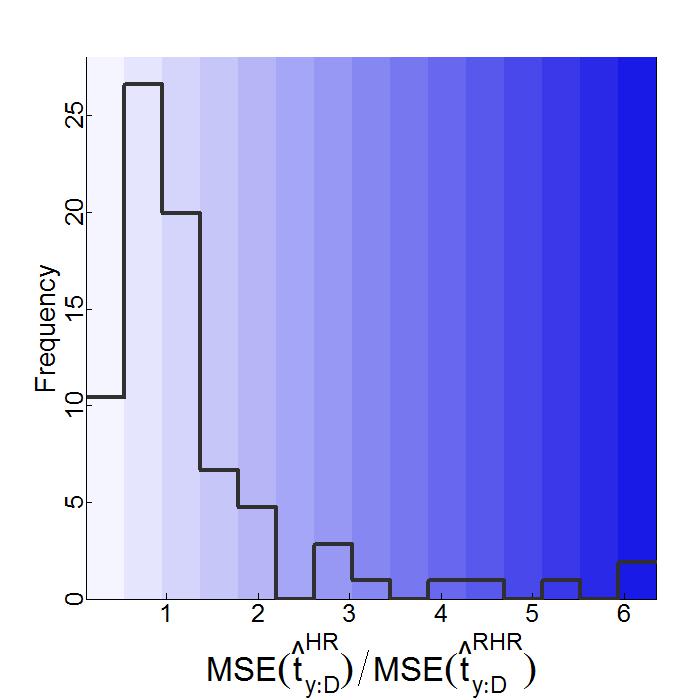}
  \label{fig:p3c_1}
\end{subfigure}%
\begin{subfigure}{.25\textwidth}
  \centering
  \includegraphics[width=0.9\linewidth]{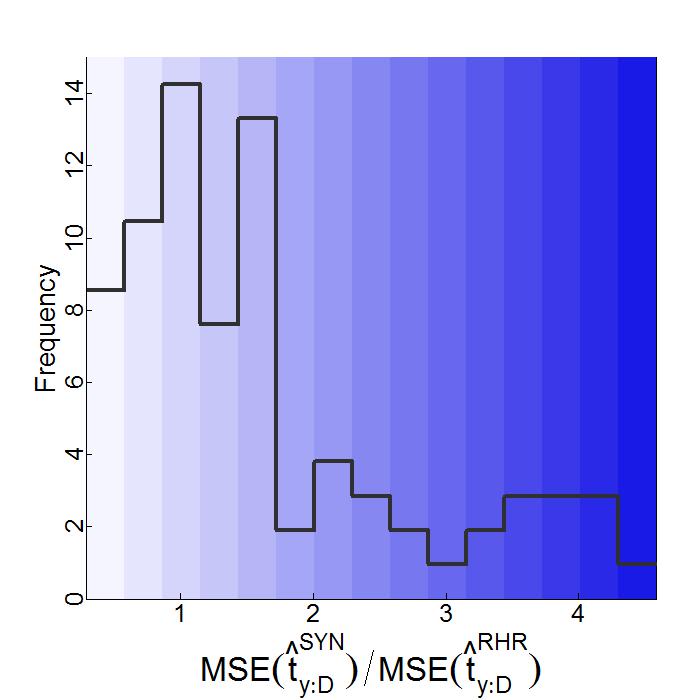}
  \label{fig:p3c_2}
\end{subfigure}%
\begin{subfigure}{.25\textwidth}
  \centering
  \includegraphics[width=0.9\linewidth]{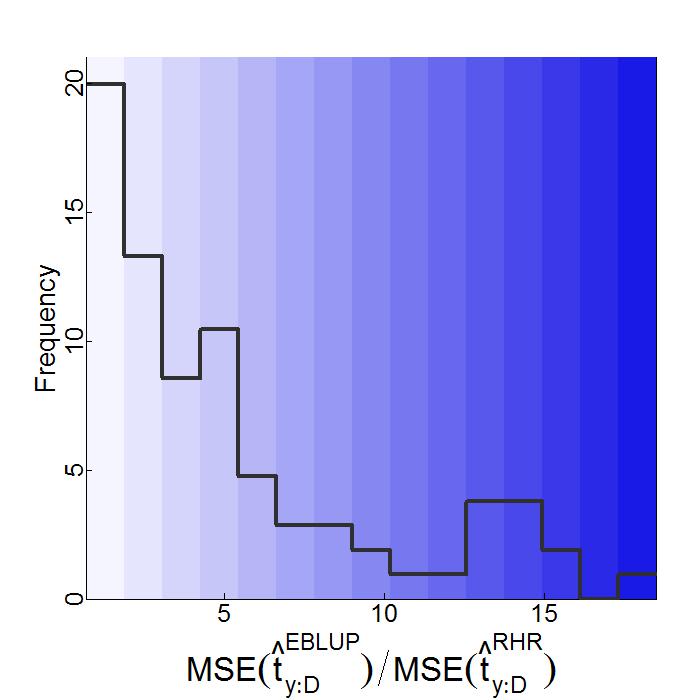}
  \label{fig:p3c_3}
\end{subfigure}%
\begin{subfigure}{.25\textwidth}
  \centering
  \includegraphics[width=0.9\linewidth]{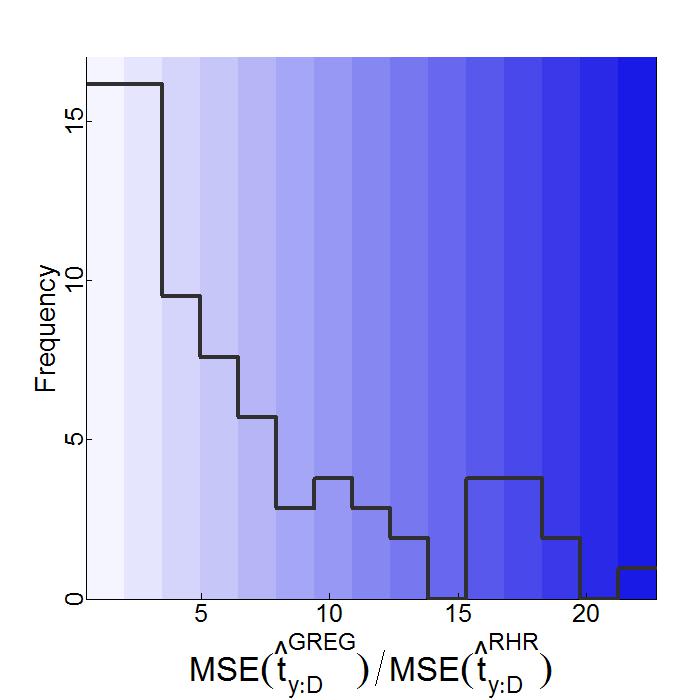}
  \label{fig:p3c_4}
\end{subfigure}
\caption{Counts for ($P_2$), Case (C), and HR, SYN, EBLUP, GREG vs RHR}
\label{fig:p3c_hist}
\end{figure}

\clearpage{}
\newpage{}

\small

\bibliographystyle{plain} 
\bibliography{liter_abc}
 
\end{document}